\documentclass[a4paper,10pt]{article}
\pdfoutput=1

\usepackage{amsmath}
\usepackage{amssymb}

\usepackage{placeins}

\usepackage[svgnames]{xcolor}

\usepackage{multicol}
\usepackage{makeidx}
\makeindex

\usepackage{microtype}

\newcommand{\jjref}[1]{\ref{#1} on page \pageref{#1}}%

\newcommand{\jjterm}[1]{\emph{#1}\index{#1}\marginpar{\small\emph{#1}}}
\newcommand{\jjxterm}[2]{\emph{#1}\index{#2}\marginpar{\small\emph{#2}}}
\newcommand{\eqq}[1]{``{#1}''}
\newcommand{\Lmap}[2]{\texttt{#1} $\mapsto$ \texttt{#2}}

\newcommand{\adeg}[0]{^{\circ}}

\def\jjstab{\mathop{\mathrm{Stab}}\nolimits}%
\def\jjorb{\mathop{\mathrm{Orb}}\nolimits}%

\newcommand{\tileS}[2]{[\texttt{#1}]$^{#2}$}

\makeatletter
\renewcommand{\fps@figure}{h!tbp}
\renewcommand{\fps@table}{h!tbp}
\makeatother

\usepackage{parskip}
\usepackage{subcaption}

\usepackage{graphicx}
\graphicspath{ {./img-pdf/} }

\newcommand{\ifDEBUG}[1]{}

\newcommand{\IMAGE}[1]%
{\input{figures/#1.tex}}

\ifDEBUG
{
\renewcommand{\IMAGE}[1]%
{\input{figures/#1.tex}}
}
\makeatletter
\renewcommand{\ps@headings}{%
\renewcommand{\@evenfoot}{}%
\renewcommand{\@oddfoot}{}%
\renewcommand{\@oddhead}{}%
\renewcommand{\@evenhead}{}%
}%
\let\ps@plain=\ps@headings 
\makeatother
\usepackage{fancyhdr}
\usepackage{ifpdf}
\ifpdf
\usepackage[hyperindex,pdftex,breaklinks=true,backref=page,colorlinks=true,linkcolor=DarkBlue,citecolor=DarkBlue,urlcolor=DarkBlue]{hyperref}%
\hypersetup{%
pdfauthor={Joerg Arndt, Julia Handl},%
pdftitle={Edge-covering plane-filling curves on grid colorings},%
pdfsubject={plane-filling curves},%
pdfkeywords={plane-filling curves},%
bookmarksdepth=subsubsection,%
bookmarksopen,
bookmarksopenlevel={3},
}
\else
\usepackage[hypertex,hyperindex,breaklinks=true,backref=page,colorlinks=true,linkcolor=DarkBlue,citecolor=DarkBlue,urlcolor=DarkBlue]{hyperref}%
\fi

\makeatletter
\@addtoreset{equation}{section}

\@addtoreset{figure}{subsubsection}

\makeatother

\begin{document}
\bibliographystyle{plain}
\title{Edge-covering plane-filling curves on grid colorings:
\\ a pedestrian approach
}
\author{J\"{o}rg Arndt, Julia Handl}
%
%
%
%
\maketitle

\typeout{jjstart}

\begin{abstract}\noindent
We describe families of plane-filling curves on
any edge-to-edge tiling of the plane with regular polygons
and finitely many classes of edges.
It is shown how to partition the minimal
number of edge classes
from the group $G$ of symmetries of the tiling
into refined colorings of the tiling,
corresponding to finite subgroups of $G$.
All of these colorings correspond to
families of plane-filling curves which
we call curve-sets.
Our exposition is driven by illustrated examples.

%
\end{abstract}

\clearpage
{
\tableofcontents
}

\clearpage

\section{Terminology and definitions}
We fix some terminology.

\subsection{\texorpdfstring{$k$}{k}-uniform grids and their edge classes}


\IMAGE{easy-grids}

An edge-to-edge tiling (also called tessellation or grid)
by regular polygons is commonly
defined by specifying the classes of points (vertices).
Four grids with just one class of points are shown in Figure~\ref{fig:easy-grids},
the respective symbols $(4^4)$, $(3^6)$, $(3.6.3.6)$, and $(3.4.6.4)$
specify the numbers of edges of the polygons adjacent to each point.
Among the possible cyclic rotations the lexicographic minimum is chosen.

A grid with two point classes is shown in
Figure~\ref{fig:3446-3464-grid-graph}.
%
The symmetries of the grid usually allowed
are translations, rotations, and reflections.
A tiling with $k$ classes of points is called $k$-uniform,
see \cite[Section~2.2, pp.~65ff]{gruenbaum-shepard-1986}.

\IMAGE{3446-3464-grid-graph}

For our purposes, we are rather interested in
the \jjxterm{edge classes}{edge class}.
Moreover, we disallow reflections and use directed edges.
We use one color and the equivalent letter for each edge class
in the grid drawing.

We use the term \jjterm{grid} also
to refer to the \emph{directed} grids from here on.

A grid can be defined using the \jjterm{transitions}
between successive edges.
The transitions are triples $(F, t, G)$
with (edge-)letters $F$ and $G$ and turn $t$ have
the property that the pair $(F, t)$ uniquely defines $G$.
We call this the \jjterm{unique transition property}.
In other words, the first two entries in a transition determine the third.
It is also true that the last two entries determine the first.
In our $(3.4.4.6; 3.4.6.4)$-grid shown in Figure~\ref{fig:3446-3464-grid-graph}
the transitions can be given as
\begin{verbatim}
  A++A,  B++++C, C++++D, D++++B, E+++F, F+++E
  A---B, B---F,  C--E,   D---A,  E--C,  F---D
\end{verbatim}
where a word of $j$ letters \texttt{+} or \texttt{-} respectively is a
counterclockwise or clockwise turn by $j\,2\,\pi/12=$ $j\cdot{}30$ degrees.
We will use the abbreviations
\jjxterm{CW}{CW sense of rotation} and
\jjxterm{CCW}{CCW sense of rotation}
respectively for clockwise and counterclockwise
sense of rotation from this point on.

Another way of defining a grid is
by specifying the \jjxterm{prototiles}{prototile}.
The prototiles of the $(3.4.4.6; 3.4.6.4)$-grid
are shown on the right in Figure~\ref{fig:3446-3464-grid-graph}.
We denote the CCW tiles of this grid
by
 \tileS{A++}{6},
 \tileS{B++++C++++D++++}{1}, and
 \tileS{E+++F+++}{2}
and the CW tiles by
 \tileS{A---B---F---D---}{1} and
 \tileS{C--E--}{3}.
The exponent at the right indicates a repetition,
for example, \tileS{E+++F+++}{2} could also
be written as \tileS{E+++F+++E+++F+++}{}.
%
Note that a prototile can appear in more then one orientation.
For example, note how the prototile \tileS{A---B---F---D---}{1}
(top right in Figure~\ref{fig:3446-3464-grid-graph})
appears in six orientations around the blue hexagons.

A \jjterm{lattice tile} is a tile that tiles the plane
in such a way that it appears in just one orientation.
For brevity we sometimes just say \jjterm{tile}
instead of prototile, (whether or not it is a lattice tile).
This will never lead to confusion, if a prototile
in a grid coloring is a lattice tile, we say so.

\IMAGE{easy-grids-dir}

If there is just one CW (or CCW) tile and it appears in just
one orientation, then it is a lattice tile:
it can be decomposed into smaller copies of
itself, each of the same orientation;
repeating this process shows that it is indeed a lattice tile.
This happens for both tiles of the square $(4^4)$ grid and
the triangle $(3^6)$ grid,
and for the hexagonal tile of the trihexagonal grid $(3.6.3.6)$,
see Figure~\ref{fig:easy-grids-dir}.
The triangular tile of the trihexagonal grid
appears in two orientations, therefore is not a lattice tile.

We'll give examples where certain arrangements of
several prototiles give a lattice tile.



\subsection{Edge-covering curves and their tiles}

For the description of curves on a grid we use Lindenmayer systems
as described in
\cite[Section~1.2, pp.~3ff]{arndt-curve-search}.
%
A \jjterm{Lindenmayer system} (or L-system)
is a triple $(\Omega, A, P)$ where $\Omega$ is an alphabet,
$A$ a word over $\Omega$ (called the \jjterm{axiom}),
and $P$ a set of maps from letters $\in{}\Omega$
to words over $\Omega$ that contains one map for each letter.

The word that a letter is mapped to is called the
\jjterm{production} of the letter.
We also use the term \jjterm{map} for this.
If the production of a letter is just the letter itself,
we call the letter a \jjterm{constant}.

We specify curves by L-systems interpreted as a sequence of
unit-length edges and turns.
The curves can be rendered via turtle graphics,
see \cite[Section~1.3, pp.~6ff]{algorithmic-beauty-of-plants}.
%
The initial direction is arbitrary.
Letters are interpreted as \eqq{draw a unit stroke in the current direction},
\texttt{+} and \texttt{-} as turns by $\pm$ a fixed angle $\varphi=\frac{2\,\pi}{n}$
where $n \in \mathbb{N}_+$.
We will also use the letter \texttt{0} for turns by $0\adeg$ (non-turns).


Applying the maps of the L-system $j$ times to the axiom
gives a word that is the $j$th iterate.
We call the drawing corresponding to this word the $j$th \jjterm{iterate} of the curve.
The first, second, and fourth iterates for the axiom \texttt{F} (corresponding to a single edge)
with production \Lmap{F}{F0F+F0F0F-F+F-F-F0F+F+F-F} are given in
Figure~\ref{fig:333333-curve-r13-1}.
The first iterate is also called the \jjterm{motif},
the motif for a letter $F$
corresponds to the drawing of the production of $F$.
The letters \texttt{+} and \texttt{-} respectively
denote turns by $\pm{}120$ degrees,
and \texttt{0} for a (non-)turn by 0 degrees.
The curve is drawn with rounded corners to make the sequence
of strokes apparent.

\IMAGE{333333-curve-r13-1}

Some sources, for example
\cite{mandelbrot-fractal-geom-nature}
and
\cite{peitgen-chaos-and-fractals},
call the axiom the initiator
and use the term generator for the (single) map of
a Lindenmayer system with only one non-constant symbol.


The colors used in Figure~\ref{fig:333333-curve-r13-1}
do \emph{not} correspond to edge classes;
the curve lives on the triangle grid where there is just one class.
Instead, each stroke of the motif has its own (arbitrary) color
and the parts of the higher iterates of the curve inherit the
color of the edge it developed from.
This makes the \jjterm{self-similarity} of the curve obvious,
each colored part in the drawing on the right has the same shape,
and magnifying those shapes gives the shape of the entire drawing.

All of our curves will be \jjterm{self-avoiding}
in the sense that both self-crossings
and doubly drawn edges are forbidden.

We call a (self-avoiding) curve \jjterm{edge-covering}
if a large enough iterate traverses all edges in
an arbitrarily large disc on the grid.
Such a curve is also \jjterm{plane-filling}.
There are plane-filling curves of other types as well, for example,
curves that traverse all points in such a disc,
well-known examples being the Hilbert curve and the Peano curve.

\IMAGE{333333-curve-r13-1-solid}
Figure~\ref{fig:333333-curve-r13-1-solid}
shows a rendering of the curve
already seen in Figure~\ref{fig:333333-curve-r13-1}
where each edge is drawn as a lozenge
stretching to the midpoint of the neighboring triangles.
We call this an \jjterm{area drawing} (for the lack of a better term).
Area drawings are useful for emphasizing the shape of a curve and we will
use it for many images without saying so.  Indeed, the fourth iterate
on the right in Figure~\ref{fig:333333-curve-r13-1} is an area drawing.
The idea for drawing curves this way appeared
in \cite[Addendum, pp.603-605]{davis-knuth-new}.

The number of edges incident to a point (in either direction)
will be called the \jjterm{valency} of the point.
It appears that each point needs to have an even valency,
because for every stroke to the point, a stroke away from the point is needed.
We will, however, solve this apparent problem
in Section~\jjref{sect:double-edges}.


\FloatBarrier

\subsection{Curves on grids with more than one edge class: curve-sets}\label{sect:curve-sets}

Edge-covering curves on grids with just one edge class
have been treated in~\cite{arndt-curve-search}.
For grids with $j>1$ edge classes we have $j$ curves working together to fill the plane,
we call this a \jjterm{curve-set}.



\IMAGE{3464-tile}

Our first example of a curve-set lives on
the $(3.4.6.4)$-grid with $j=2$ edge classes,
shown in Figure~\ref{fig:3464-tile}.
The classes are denoted by the letters
\texttt{A} (edges surrounding hexagons, drawn in blue)
and \texttt{B} (edges surrounding triangles, drawn in red).
The three prototiles are given by the words
\tileS{A++}{6} = \texttt{A++A++A++A++A++A++} (hexagon),
\tileS{B++++}{3} = \texttt{B++++B++++B++++} (triangle), and
\tileS{A---B---}{2} = \texttt{A---B---A-{}--B---} (square).
The turns corresponding to \texttt{+} and \texttt{-}
are by $\pm{}30$ degrees ($2\,\pi/12$).
The transitions of the grid are
\texttt{A++A}, \texttt{B++++B}, \texttt{A---B}, and \texttt{B---A},
as can be read from the prototiles.

The two curves of the curve-set have the maps
{\small
\begin{verbatim}
  A |--> A++A++A++A---B---A---B++++B++++B---A++A++A++A--- \
          B---A---B++++B---A---B++++B---A---B++++B++++B---A
  B |--> B---A++A++A---B++++B++++B---A---B---A++A---B++++B
\end{verbatim}
}
%
The first iterate of both curves are shown in
Figure~\ref{fig:3464-ju-r19-a-curves-it1}.


The curves of a curve-set
exhibit \jjterm{mutual self-similarity}
as shown in Figure~\ref{fig:3464-ju-r19-a-curves-decomp}:
both curves are decomposed into smaller copies of each other.
The decomposition into colors is done similarly as in
Figure~\jjref{fig:333333-curve-r13-1}.

\IMAGE{3464-ju-r19-a-curves-it1}
\IMAGE{3464-ju-r19-a-curves-decomp}

The \jjterm{substitution matrix} of a curve-set is the matrix
whose entry $r,c$ is the number of letters $r$ in the production of class $c$.
For this curve-set the matrix is
\[
\left[ {\begin{array}{cc}
    13 & 6 \\
    12 & 7 \\
\end{array} } \right]
\]
It tells us that curve \texttt{A}
(left in Figure~\ref{fig:3464-ju-r19-a-curves-decomp})
is decomposed into 13 smaller copies of curve \texttt{A} (itself) and
12 smaller copies of curve \texttt{B} (left column of the matrix).
Curve \texttt{B} (right in Figure~\ref{fig:3464-ju-r19-a-curves-decomp})
is decomposed into 6 copies of curve \texttt{A}
and 7 copies of curve \texttt{B} (itself).
This can be seen in Figure~\ref{fig:3464-ju-r19-a-curves-it1} as well.

We observe that all row sums $R$ of the matrix are equal ($R=19$).
This allows us to define the \jjterm{order} of a curve-set to be $R$.
The definition of the order in \cite[bottom of p.4]{arndt-curve-search}
is a special case of this definition.


\subsubsection{All edge classes in a patch have the same cardinality}
We show that the row sums of the substitution matrix
are indeed always identical,
hence our concept of order is well-defined.

\IMAGE{zout-2x3-max-cols-ed}
Consider a parallelogram (minimal or not) satisfying the
translational symmetries of the grid.
We call such a parallelogram composed of $r\times{}s$ minimal parallelograms
an $(r,s)$--\jjterm{patch}.
Figure~\ref{fig:zout-2x3-max-cols-ed} shows $(2,3)$--patches on the square grid.
A patch has the same rotational symmetries as the grid.
The group of the translational symmetries of the patch is finite;
it is the direct product of two cyclic groups of orders $r$ and $s$.

Let $G$ be the group of all symmetries of the edges for a patch.
Consider the group action $G\times{}X\to{}X$ where $X$ is
the set of the edges in the patch.
Observe that no element in $G$ apart from the identity fixes any edge.
This is trivially true for translations, but also true for rotations
as we consider \emph{directed} edges and disallow reflections.
The set of group elements fixing an element (edge) $x$ is called the
stabilizer of $x$, denoted as $\jjstab(x)$.
We just observed that the stabilizer of all elements contains just one element,
the neutral element of the group, so $\left|\jjstab(x)\right|=1$.

%
Burnside's lemma (see \cite[Theorem~25.17, p.~293-294]{lidl-pilz-aaa})
tells us that
\[
\left|\jjorb(x)\right| \cdot \left|\jjstab(x)\right|
  = \left|G\right| \qquad \forall x \in X
\]
where $\left|\jjorb(x)\right|$ is the number of edges equivalent to the edge $x$
and $\left|G\right|$ the number of elements of $G$.
This is called the Orbit-Stabilizer Theorem,
see \cite[Theorem~(17.2) and Corollary~(17.3), p.~94]{armstrong}.
Now $\left|\jjstab(x)\right|=1$, so
\[
\left|\jjorb(x)\right| = \left|G\right| \qquad \forall x \in X
\]
In other words, all equivalence classes have the same size,
namely $\left|G\right|$.
Note that our argument falls apart when
either using undirected edges or including reflections in the group.

\subsubsection{The prototiles of a grid}

The prototiles of a grid fall into two classes,
some have CCW sense of rotation
(tiles \tileS{A++}{6} and \tileS{B++++}{3} in
Figure~\jjref{fig:3464-tile}),
others have CW sense (tile \tileS{A---B---}{2} in
Figure~\ref{fig:3464-tile}).
The first iterates (with respect to our curve-set)
of the CCW prototiles are shown in
Figure~\ref{fig:3464-ju19-curve-prototiles-plus-it1}.

\IMAGE{3464-ju19-curve-prototiles-plus-it1}

\IMAGE{3464-ju19-curve-prototiles-plus-decomp}
In the limit, each prototile obtains a shape determined by the curve-set,
see Figure~\ref{fig:3464-ju19-curve-prototiles-plus-decomp}.
These shapes are mutually self-similar:
the shape of the tile \tileS{A++}{6} is that of the smaller blue regions,
the shape of \tileS{B++++}{3} is that of the red regions.
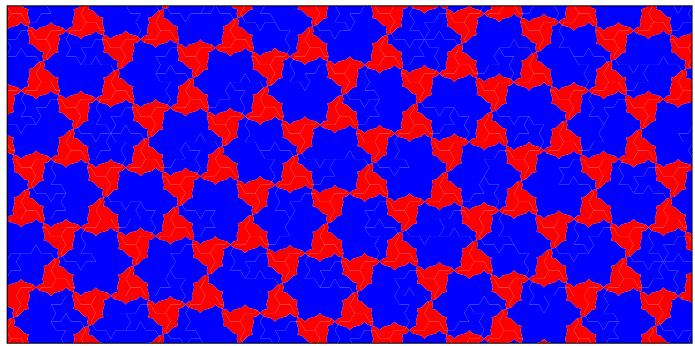
The two CCW prototiles tile the plane as shown in
Figure~\ref{fig:3464-ju19-curve-prototiles-plus-tiling}.
\IMAGE{3464-ju19-curve-prototiles-plus-tiling}

There is just one CW prototile, shown in
Figure~\ref{fig:3464-ju19-curve-prototile-minus}.
Its shape is self-similar, it can be decomposed into 19 smaller copies of itself
(colored by orientation, right in figure).
The CW prototile tiles the plane, occurring in three orientations.
\IMAGE{3464-ju19-curve-prototile-minus}

What we observed here is always true.
All prototiles of one sense of rotation (CW or CCW)
are mutually self-similar.
Also, all prototiles of one sense of rotation
together tile the plane.

The exponent in the notation for a prototile
is the guaranteed rotational symmetry of the shape of the tile.
For example, the tile \tileS{B++++}{3} will have
at least 3-fold rotational symmetry.
Sometimes greater symmetries are possible, as we will show later.
For brevity, we sometimes call a $k$-fold rotational symmetry just
\jjterm{$k$-fold symmetry}.


\FloatBarrier

\subsection{Refined grid colorings}

The grid colorings considered so far correspond to
the full group $G$ of the symmetries of the grid.
These are the \jjxterm{minimal colorings}{minimal coloring}
of the grid, containing the least number of classes.
Using a finite subgroup of $G$
one obtains \jjxterm{refined colorings}{refined coloring},
containing more classes.

One can always put each edge in a patch into its own class.
Figure~\jjref{fig:zout-2x3-max-cols-ed} shows such a maximally refined
coloring for $(2,3)$--patches on the square grid.
The maximally refined coloring
on the square grid was used in
\cite[pp.~132--134]{arndt-handl-bridges2019}
 for $(1,1)$--patches
and
\cite[pp.~136--138]{arndt-handl-bridges2019}
 for $(1,2)$--patches.

\subsubsection{A first example: a coloring of the square grid with two colors}\label{sect:first-ex-4444-lr-col}

\IMAGE{4444-rotational-lr-grid-graph}
As an example of a refined coloring we use the square grid.
We reduce the 4-fold rotational symmetry of the square grid to a 2-fold symmetry
to obtain the coloring shown in Figure~\ref{fig:4444-rotational-lr-grid-graph}.
We assign the letters \texttt{L} and \texttt{R} respectively to
the horizontal and vertical edges.
The prototiles are \tileS{L+R+}{2} and \tileS{L-R-}{2}.
The transitions for this coloring are
\texttt{L+R},
\texttt{L-R},
\texttt{R+L}, and
\texttt{R-L}.

The first iterates of a curve-set of order 13 with maps
\Lmap{L}{L-R-L+R-L+R+L+R-L+R-L-R+L} and
\Lmap{R}{R+L-R-L-R+L+R+L-R-L-R+L+R}
are shown in Figure~\ref{fig:4444-counter-example-curve-r13-x01-it1}.
The decompositions of the curves for \texttt{L} and \texttt{R}
into smaller copies of themselves are
shown in Figure~\ref{fig:4444-counter-example-curve-r13-x01-decomp}.
\IMAGE{4444-counter-example-curve-r13-x01-it1}
\IMAGE{4444-counter-example-curve-r13-x01-decomp}

\FloatBarrier
Dekking treated \jjterm{folding curves} on this coloring
of the square grid \cite{dekking-tiles-final}.
For those curve-sets the production of \texttt{R} is obtained from
the production of \texttt{L} by reversing it and swapping
turns \texttt{+} with \texttt{-}
and
letters \texttt{L} with \texttt{R}.
In other words, the curve for \texttt{R} is that of \texttt{L},
in reversed order.
The curve we just presented is not a folding curve.

\IMAGE{folding-r09-it1}
The first iterates of a folding curve of order 9 with maps
\Lmap{L}{L+R+L-R+L+R-L+R-L}
and
\Lmap{R}{R+L-R+L-R-L+R-L-R}
are shown in Figure~\ref{fig:folding-r09-it1}.
The shapes of either curve (both \texttt{L} or \texttt{R})
and either tile (\tileS{L+R+}{2} or \tileS{L-R-}{2})
are shown in Figure~\ref{fig:folding-r09-shapes}.
That the shapes of both curves coincide
is a consequence of being folding curves.
%
\IMAGE{folding-r09-shapes}



\IMAGE{folding-r09-lattice-tile}
There is just one CCW prototile and
it appears in just one orientation, so it is a lattice tile.
This is shown in Figure~\ref{fig:folding-r09-lattice-tile}.
The same is true for the CW tile, which (in the special case
of folding curves!) is identical.

\subsubsection{Dekking's tile criteria}\label{sect:dekking}

In \cite{dekking-tiles-final}, Dekking gives two theorems about folding curves.

\paragraph{Dekking-1:}
If the first iterate of both tiles (\tileS{L+R+}{2} and \tileS{L-R-}{2})
are self-avoiding, then all iterates are self-avoiding.
In other words, the folding curve is self-avoiding.

\IMAGE{folding-r09-non-filling}
As an example, we use the self-avoiding (but not plane-filling)
curves with maps
\Lmap{L}{L-R+L+R-L} and
\Lmap{R}{R+L-R-L+R}.
The first and third iterates of the tile \tileS{L+R+}{2}
are shown in Figure~\ref{fig:folding-r09-non-filling}.

\paragraph{Dekking-2:}
If the interiors of both tiles are filled, the curve is plane-filling.

\IMAGE{folding-r09-plane-filling}
The tile \tileS{L+R+}{2} for the plane-filling curves with maps
\Lmap{L}{L+R-L-R-L+R+L+R-L} and
\Lmap{R}{R+L-R-L-R+L+R+L-R}
is shown in Figure~\ref{fig:folding-r09-plane-filling}.


The following generalization of Dekking's tile criteria for curve-sets
appears natural.

\paragraph{Dekking-1 CS:}
If the first iterates of all prototiles of a grid are self-avoiding
then all iterates are self-avoiding.
In other words, the curve-set is self-avoiding.

The proof Dekking gave can be extended to the
more general situation of curve sets
[Dekking, 30-November-2023, personal communication].

Note that the condition is equivalent to demanding
that the first iterates of all transitions are self-avoiding.


\IMAGE{sausage-problem}
The following is \emph{wrong} for curve-sets.
\vspace*{-1em}
\paragraph{Dekking-2 CS (A):}
If the interiors of all prototiles are filled,
the curve-set is plane-filling.

Dekking \cite[Figure~15, p.~14]{dekking-1981}
gave an example where the second criterion does not work.
We give a slightly simpler one.
Consider the curve with maps
\Lmap{L}{L+R-L} and \Lmap{R}{R}.
The first and third iterates of the prototile \tileS{L+R+}{2}
are shown in Figure~\ref{fig:sausage-problem}
(iterates of the other prototile \tileS{L-R-}{2}
give essentially identical images).
The iterates of \texttt{L} give arbitrarily long zigzag lines
(and \texttt{L} is constant).
This is certainly not a plane-filling curve.

One would hope for the following to be true.
\vspace*{-1em}
\paragraph{Dekking-2 CS (B):}
If no letter is constant and
the interiors of all prototiles are filled,
the curve-set is plane-filling.

Sadly, this is \emph{still} wrong.
A counterexample is the curve-set
(on a coloring of the square grid treated in
Section~\jjref{sect:4444-set-coloring})
with maps
\begin{verbatim}
  A |--> A+B-A
  B |--> B-A+B-A+B-A+B-A+B
  C |--> C+D-C+D-C+D-C+D-C
  D |--> D-C+D
\end{verbatim}

\IMAGE{dekking-2-omg}
As shown in Figure~\ref{fig:dekking-2-omg},
the situation is essentially the same as in the last counterexample:
all four curves and the two tiles have the limit shape of a line.
This curve-set was found by Daniel Fischer~\cite{dan-master-thesis}.

The situation is not dire, though.
Curve-sets that are not plane-filling
are recognized easily.


The substitution matrices of
both (counter)examples are reducible
(as defined in
 \cite[Definition~6.2.21, p.~403]{horn-johnson-matrix-analysis}).
For the last example it is
\[
\left[ {\begin{array}{cccc}
 2 & 4 & 0 & 0 \\
 1 & 5 & 0 & 0 \\
 0 & 0 & 5 & 1 \\
 0 & 0 & 4 & 2 \\
\end{array} } \right]
\]
This is clearly reducible,
it even is of block-diagonal form.
For the example with one letter constant the
substitution matrix is
reducible as well.
%

It appears that an irreducible
substitution matrix implies
that a curve-set is plane-filling.
However, this is not a necessary condition:
we'll present plane-filling curve-sets with
one or more letters constant;
their substitution matrices are reducible.


\FloatBarrier

\subsubsection{Creating refined grid colorings}

For the sake of simplicity, we use the square grid for our description
of a systematic way of finding colorings.

For each coloring there is a minimal vector
$\vec{s}=[r,\,c]$ with $r,\,c > 0$
such that every square is moved to a copy of itself
by any shift (along diagonals)
 $[i\cdot{}r,\, j\cdot{}c]$ with $i,\,j\in{}\mathbb{Z}$.
Figure~\jjref{fig:zout-2x3-max-cols-ed}
shows a coloring with 24 colors
where $\vec{s}=[2,\,3]$.

\subsubsection*{A computer search}

A systematic search for colorings can proceed as follows.
Choose a vector $\vec{s}=[R,\,C]$ with $R,\,C > 0$
and work on the toroidal graph obtained by connecting
edges at opposite sides as required.
Tentatively assign colors to edges
in all ways satisfying the unique transition property.
Some of the colorings found are duplicates of another coloring
after a suitable permutation of letters; these must be discarded.

The search will find all colorings whose minimal vectors $\vec{s}=[r,\,c]$
are such that $r$ divides $R$ and $c$ divides $C$.

Let $G$ be the group of symmetries of the toroidal graph.

One will always find the maximally refined coloring with all edges
in different classes, corresponding to the trivial subgroup of $G$
consisting of just the neutral element
(see Figure~\jjref{fig:zout-2x3-max-cols-ed}).
For this coloring we have $\vec{s}=[R,\,C]$.

One will also always find the minimal coloring corresponding to
the whole group $G$ (for the square grid just one edge class).
For this coloring we have $\vec{s}=[1,\,1]$.

\IMAGE{grid-colorings-1}
We did such a search for the square grid.
There are two colorings using two colors,
as shown in Figure~\ref{fig:grid-colorings-1}.
The vector $\vec{s}$ is given for each coloring.

\IMAGE{grid-colorings-2}
\IMAGE{grid-colorings-3}
\IMAGE{grid-colorings-4}
Figures~\ref{fig:grid-colorings-2},
\ref{fig:grid-colorings-3}, and
\ref{fig:grid-colorings-4}
show all five colorings with four colors
and the one coloring with five colors.

\subsubsection*{Ad hoc methods}

Grid coloring with many colors tend to lead to
curve-sets that look just random,
see \cite[Figure~15, p.~137]{arndt-handl-bridges2019}
for one example on the square grid using eight colors.

So one may focus on colorings with only a few colors.
Some ad hoc constructions tend to work.
For example, the coloring of the triangle grid
shown in Figure~\jjref{fig:333333-rotational-grid-graph}
was obtained by dropping the rotational symmetry
around the center of the triangles.
To verify a tentative coloring, one just has
to check the unique transition property.

Indeed, we obtained all grid colorings in an ad hoc manner,
except for the square grid.


\FloatBarrier
\clearpage

\section{Examples of curve-sets}

A relentless onslaught of examples lies ahead of you.
Enjoy.



\subsection{The square grid}

Here we use turns by $90$ degrees ($2\,\pi/4$).

\subsubsection{The minimal coloring with one color}


\IMAGE{4444-r29-84-dragon}

Examples of curves on the minimal coloring
are given in \cite[Section~6.2, pp.60--64]{arndt-curve-search}.
A curve with 2-fold rotational symmetry is shown in
Figure~\ref{fig:4444-r29-84-dragon}.
The coloring is by the edges of the motif in the left.
The map of the curve is
\begin{verbatim}
  F |--> F+F+F-F-F-F+F+F-F-F-F+F+F-F-F+F+F-F-F+F+F+F-F-F+F+F+F-F-F
\end{verbatim}

\subsubsection{A coloring with two colors}

\IMAGE{4444-translational-grid-graph}

We use the coloring shown in Figure~\ref{fig:4444-translational-grid-graph}
(this is the coloring shown on the left in Figure~\jjref{fig:grid-colorings-1}).
Our curve-set of order 25 has the following maps:
\begin{verbatim}
  F |--> F+F+F-G-F+F-G-F-G+G+G-F+F-G-F-G+G+G+G-F+F+F-G+G-F
  G |--> G-F+F-G+G+G-F+F+F-G-F+F+F-G-F-G+G+G-F-G-F+F-G+G+G
\end{verbatim}
\IMAGE{4444-translational-motif}
\IMAGE{4444-translational-decomp}

The first iterates of the two curves are shown in Figure~\ref{fig:4444-translational-motif}.
These curves are mutually self-similar, see Figure~\ref{fig:4444-translational-decomp}.

\IMAGE{4444-translational-tile-FGFG-self-sim}
The CW prototile \tileS{F-G-}{2} is self-similar, it can be decomposed into
smaller copies of itself, appearing in two orientations,
as shown in Figure~\ref{fig:4444-translational-tile-FGFG-self-sim}.



\IMAGE{4444-translational-r25-3s-tiles-plus}
The CCW prototiles are \tileS{F+}{4} and \tileS{G+}{4}.
Their iterates will have 4-fold rotational symmetries and be mutually self-similar.
However, for this curve-set the iterates of the two CCW prototiles are squares,
so the mutual self-similarity is boring (the $5\times{}5$ chessboard pattern).
For an example of non-trivial mutual self-similarity
see Figure~\ref{fig:4444-translational-r25-3s-tiles-plus}.
The curve-set is of order 25 and has maps
{\small
\begin{verbatim}
  F |--> F+F+F-G-F+F-G-F-G+G+G-F+F-G-F-G+G+G+G-F+F+F-G+G-F
  G |--> G-F+F-G+G+G-F+F+F-G-F+F+F-G-F-G+G+G-F-G-F+F-G+G+G
\end{verbatim}
}

\IMAGE{4444-r49-double-manta-it1}
Curve sets where both curves have reflection symmetry exist,
an example with order 49 is shown in
Figures~\ref{fig:4444-r49-double-manta-it1} and \ref{fig:4444-r49-double-manta}.
The maps are
{\scriptsize
\begin{verbatim}
  F |--> F+G-F+G+F-G-F-G+F-G-F+G-F+G+F-G-F+G+F+G-F+G+F-G-F+G-F-G-F+G+F+G-F+G-F-G-F+G+F+G-F
  G |--> G+F-G+F+G-F-G+F-G-F-G+F+G+F-G-F-G+F-G+F-G+F+G-F-G+F+G-F+G+F+G-F-G-F+G+F+G-F+G-F-G+ \
         F-G-F-G+F+G+F-G+F-G-F-G+F+G+F-G
\end{verbatim}
}
\IMAGE{4444-r49-double-manta}

\FloatBarrier
\subsubsection{The other coloring with two colors}

This is already treated in Section~\jjref{sect:first-ex-4444-lr-col},
here we just give a numeration system corresponding to the tile
of a folding curve.

\subsubsection*{Complex numeration system}\label{sect:numsys}

The numeration system for the folding curve with maps
\begin{verbatim}
  L |--> L+R+L-R+L+R-L+R-L
  R |--> R+L-R+L-R-L+R-L-R
\end{verbatim}
whose motifs are shown in Figure~\jjref{fig:folding-r09-it1}.

\IMAGE{folding-r09-tile-numsys}
The first iterate of the tile is shown on the left in
Figure~\ref{fig:folding-r09-tile-numsys}.
We take the central square as the origin of the complex plane.
The coordinates of all 9 squares give the set
\[D = \left\{
 0,\;
+1-i,\; -1+i,\;
+2 i,\; -2 i,\;
+1+3 i,\; -1-3 i,\;
+2+2 i,\; -2-2 i
 \right\}\]
A numeration system with radix $+3$
and digits $D$ has the fundamental region
shown on the right.
The \jjterm{fundamental region} of a numeration system
is the set of numbers whose expansion has zero integral part
in their expansion (in that numeration system),
see \cite[Section~3.3, pp.~20-26]{arndt-curve-search}.

A complex numeration system corresponding to the tile
of a folding curve of order 5 is shown in
\cite[Figure~7, p.182]{arndt-handl-bridges2018}.

\FloatBarrier
\subsubsection{A coloring with four colors}

\IMAGE{4444-star-grid-graph}

%
With the coloring shown in Figure~\ref{fig:4444-star-grid-graph}
(same as left in Figure~\jjref{fig:grid-colorings-4})
the sequence of letters in any curve is always
\texttt{A}, \texttt{B}, \texttt{C}, \texttt{D},
\texttt{A}, \texttt{B}, \texttt{C}, \texttt{D}, \ldots,
regardless of the turns.
Each prototile appears in two orientations.

\subsubsection*{A curve-set of order 49}

%


%
Our curve-set of order 49 has the maps
{\small
\begin{verbatim}
  A |--> A+B-C+D+A-B-C-D+A+B-C+D-A-B-C+D-A+B-C+D-A+B+C-D-A+B+C+D-A-B+ \
           C-D+A+B+C-D+A-B+C-D-A+B-C-D-A+B+C+D-A+B-C-D-A+B+C+D-A
  B |--> B+C-D+A+B-C-D-A+B-C-D+A-B+C+D-A-B+C+D+A-B+C+D-A-B+C-D-A-B+C+ \
           D+A-B+C+D-A-B-C+D-A-B+C-D+A+B+C-D+A-B
  C |--> C+D-A+B-C+D-A-B-C+D-A+B-C+D+A-B+C-D+A+B+C-D+A-B-C+D-A-B-C+D+ \
           A+B-C+D-A-B-C+D+A+B-C
  D |--> D+A-B+C+D-A-B-C+D+A-B-C-D+A-B+C-D+A+B-C-D+A+B+C-D-A+B+C+D-A+ \
           B-C-D+A-B-C-D+A+B+C-D+A-B-C-D+A+B+C-D
\end{verbatim}
}
The mutual self-similarity of the four curves
is shown in Figure~\ref{fig:4444-star-r49-manta-curve-decomp}.
All curves have reflection symmetry.
\IMAGE{4444-star-r49-manta-curve-decomp}


\FloatBarrier
\subsubsection*{A curve-set of order 5 and a numeration system}

\IMAGE{4444-star-r05-lattice-tile}

Our curve-set has the maps
\begin{verbatim}
  A |--> A-B+C
  B |--> D+A-B-C-D+A-B+C+D+A-B
  C |--> C-D+A
  D |--> B-C+D
\end{verbatim}
Their motifs can be gleaned from the left depiction in
Figure~\ref{fig:4444-star-r05-lattice-tile}.
A lattice tile is created by combining two CW tiles \tileS{A-B-C-D-}{1}
into the axiom \texttt{A-B-C-D+A-B-C-D}.
The fifth iterate of this is shown in the middle of the figure.
The corresponding numeration system is shown on the right.
It has radix $-2+i$
and the digit set is
\[D = \left\{ 0,\;
 +1,\; -1,\;
 +1-i,\; -1+i
 \right\}\]
The digits correspond to the positions of the centers
of the small copies of lattice tiles.

\IMAGE{4444-star-r05-lattice-tile-plus}
We remark that the CCW lattice tile
cannot be decomposed into smaller copies of itself.
The left in Figure~\ref{fig:4444-star-r05-lattice-tile-plus}
shows an attempt to do so, starting from the two central squares.
Both the light green and the dark green parts are disconnected.
It is still a lattice tile, the whole place
can be covered with translated copies of it.
Unsurprisingly, we cannot give a
numeration system corresponding to it.

\FloatBarrier
\subsubsection{Another coloring with four colors}\label{sect:4444-set-coloring}

\IMAGE{4444-set-grid}
The coloring together with its prototiles is shown in
Figure~\ref{fig:4444-set-grid},
it is the same as shown in Figure~\jjref{fig:grid-colorings-3} (right).
Either prototile occurs in just one orientation, both are lattice tiles.

All curve-sets presented here where found
by Daniel Fischer \cite{dan-master-thesis}.

\subsubsection*{A curve-set of order 4}

\IMAGE{4444-set-dan-twindragon-square}
We consider the curve-set of order 4 with the following maps
\begin{verbatim}
  A |--> D+A-D-C-B+C+D
  B |--> A-D+A+B-A
  C |--> B+C-B
  D |--> C
\end{verbatim}
An image combining both prototiles is shown
on the right in
Figure~\jjref{fig:4444-set-dan-twindragon-square}.
The shape of the CW prototile is the twindragon,
that of the CCW prototile a square.

\subsubsection*{A curve-set of order 7}

\IMAGE{4444-set-dan-gosper-island-oid}
A curve-set of order 7 with maps
\begin{verbatim}
  A |--> A-D+A
  B |--> B-A+B
  C |--> C-B-A+B+C-B+C+D+A-D-C-B+C
  D |--> D+A-D+A-D-C-B+C+D
\end{verbatim}
The CCW prototile shown on the left in
Figure~\ref{fig:4444-set-dan-gosper-island-oid}
looks like a deformed Gosper island,
a figure that is not expected on the square grid.
The Gosper island on the right is the CCW tile
\tileS{F+}{3}
on the triangle grid (minimal coloring)
using the curve with map
\begin{verbatim}
  F |--> F+F-F-F+F+F-F
\end{verbatim}

The next example is of the same nature.

\subsubsection*{A curve-set of order 3}

\IMAGE{4444-set-dan-terdragon-oid}
A curve-set of order 3 with maps
\begin{verbatim}
  A |--> A+B-A
  B |--> B+C-B
  C |--> C+D+A-D-C
  D |--> D
\end{verbatim}
The CCW prototile is shown on the left in
Figure~\ref{fig:4444-set-dan-terdragon-oid},
only three curves are apparent because \texttt{D} is a constant.
All three curves look like distorted terdragons:
The right image shows the CCW tile
\tileS{F+}{3} of the terdragon (again on the triangle grid).
Its map is
\begin{verbatim}
  F |--> F+F-F
\end{verbatim}

\subsubsection*{A curve-set of order 6}

\IMAGE{4444-set-dan-twin-ter-oid}
Both prototiles of a curve-set of order 6
are shown in Figure~\ref{fig:4444-set-dan-twin-ter-oid}.
The maps of this curve-set are
\begin{verbatim}
  A |--> A-D-C+D-C-B+C+D+A-D+A
  B |--> B-A+B-A-D+A+B
  C |--> C-B+C-B+C
  D |--> D
\end{verbatim}
There is no valid excuse for presenting this,
but we do like the looks.




\FloatBarrier

\subsection{The triangle grid}

Here we use turns by $120$ degrees ($2\,\pi/3$).

\subsubsection{The minimal coloring with one color}

Examples of curves on the minimal coloring
are given in \cite[Section~6.1, pp.56--59]{arndt-curve-search}.

\IMAGE{333333-R25-544}

We give an example where both tiles have reflection symmetry.
The curve of order 25 with map
\begin{verbatim}
  F |--> F0F+F+F-F-F0F-F0F+F+F-F-F0F-F0F+F+F-F0F0F+F-F+F0F
\end{verbatim}
gives the tiles shown in Figure~\ref{fig:333333-R25-544}.

\clearpage
\subsubsection{A coloring with three colors}\label{sect:col-333333-translational}

\IMAGE{333333-translational-grid-graph}
We use the grid coloring
shown in Figure~\ref{fig:333333-translational-grid-graph},
together with its prototiles.

Our curve-set is of order 25 and has maps
\begin{verbatim}
  F |--> F+F-H+H-G-F0G+G0H0F
  G |--> G-F+F-H+H+H-G+G+G-F0G-F0G-F-H+H+H0F0G
  H |--> H-G+G-F+F-H+H+H0F0G-F0G+G0H-G-F+F-H+H
\end{verbatim}
The motifs of the curves
are shown in Figure~\ref{fig:333333-translational-r25-motifs}.

\IMAGE{333333-translational-r25-motifs}

\IMAGE{333333-translational-plus-tiles-self-similarity}
The three CCW prototiles,
\tileS{F+}{3}, \tileS{G+}{3}, and \tileS{H+}{3},
have 3-fold rotational symmetry
and are mutually self-similar
as shown in Figure~\ref{fig:333333-translational-plus-tiles-self-similarity}.
%

%
The arrangement \texttt{(F+F+F)-(H+H+H)-(G+G+G)} of
the CCW prototiles tiles the plane by itself (it is a lattice tile),
it is shown in the lower right


\IMAGE{333333-translational-minus-tile-self-similarity}
The one CW prototile \tileS{F-H-G-}{1} tiles the plane,
appearing in three orientations.
Its self-similarity is shown in
Figure~\ref{fig:333333-translational-minus-tile-self-similarity} (left).
The arrangement \texttt{(F-H-G)+(G-F-H)+(H-G-F)} of
the prototile \tileS{F-H-G-}{1} in its three orientations
tiles the plane,
appearing in just one orientation (so it is a lattice tile again),
see Figure~\ref{fig:333333-translational-minus-tile-self-similarity} (right).
Neither the prototile nor the arrangement have any rotational symmetry.

\FloatBarrier
\subsubsection{Another coloring with three colors}\label{sect:col-333333-star}


\IMAGE{333333-star-grid-graph}

%
%
The grid coloring is shown in Figure~\ref{fig:333333-translational-grid-graph}.
With this coloring the sequence of letters in any curve is always
\texttt{A}, \texttt{B}, \texttt{C},
\texttt{A}, \texttt{B}, \texttt{C}, \ldots,
regardless of the turns.
%

\subsubsection*{A curve-set of order 16}

\IMAGE{333333-star-r16-var1-it1}
%
Our curve-set is of order 16 and has maps
\begin{verbatim}
  A |--> A+B-C+A-B-C0A+B0C-A-B+C-A+B+C0A0B+C-A
  B |--> B+C-A+B-C-A0B+C-A+B+C0A-B-C+A-B0C+A+B-C+A0B-C-A+B
  C |--> C0A0B0C
\end{verbatim}
The first iterates are shown in Figure~\ref{fig:333333-star-r16-var1-it1}.
%

\IMAGE{333333-star-tile-plus}
\IMAGE{333333-star-tile-minus}
Each prototile appears in three orientations,
the iterates of either tile the plane.
The CCW prototile \tileS{A+B+C+}{1}
is shown in Figure~\ref{fig:333333-star-tile-plus},
the CW prototile \tileS{A-B-C-}{1}
in Figure~\ref{fig:333333-star-tile-minus}.


\IMAGE{333333-star-lattice-tiles-outline}
Here one can easily build six lattice tiles,
see Figures~\ref{fig:333333-star-lattice-tiles-outline}
and \ref{fig:333333-star-lattice-tiles}.
Each has a 3-fold rotational symmetry.
\IMAGE{333333-star-lattice-tiles}

\FloatBarrier
\subsubsection*{Another curve-set of order 16}

\IMAGE{333333-star-r16-manta-decomp}
It is possible for all curves of a curve-set
to have a reflection symmetry.
The curves of such a curve-set of order 16
having maps
\begin{verbatim}
  A |--> A0B+C-A-B0C+A+B-C+A0B-C-A+B-C0A0B+C+A-B+C-A
  B |--> B0C0A+B0C-A-B+C-A0B0C+A+B-C+A-B
  C |--> C0A0B+C-A-B0C+A+B-C
\end{verbatim}
is shown in Figure~\ref{fig:333333-star-r16-manta-decomp}.

\IMAGE{333333-star-r16-manta-grid}

The covering of the plane with second iterates of the curves
is shown in
Figure~\ref{fig:333333-star-r16-manta-grid}.

\FloatBarrier
\subsubsection{Yet another coloring with three colors}

%
\IMAGE{333333-rotational-grid-graph}
The coloring of the grid
shown in Figure~\ref{fig:333333-rotational-grid-graph}
has just one CW and CCW prototile.
Note that either prototile appears in just one orientation,
so both are lattice tiles.


\subsubsection*{Curve-set of order 7 and a numeration system}

\IMAGE{r7-x3-set-tile}
The first curve-set is of order 7 and has the maps
\begin{verbatim}
  A |--> A+B-A+B+C-B+C-B-A+B-A
  B |--> B+C0C0C-B
  C |--> C+A0A0A-C
\end{verbatim}
The first iterate of CCW tile is shown in Figure~\ref{fig:r7-x3-set-tile}.

\FloatBarrier
Now take the lower red triangle of the right image as
origin of the complex plane.
The coordinates of all 7 triangles give the set
\[D = \left\{
 0,\;
 +\omega_6,\; -\omega_6,\;
 +\omega_3,\; -\omega_3,\;
 +2{}\omega_3,\; 1 +\omega_6
 \right\}\]
where $\omega_k = \exp(2\pi{}i/k)$ (a primitive complex $k$th root of unity).
We use $D$ as digits
and the radix $-1 + 3{}\omega_6$
for the numeration system.

\IMAGE{r7-x3-set-numsys}
%
%
The fifth iterate of the CCW tile is shown on the left in
Figure~\ref{fig:r7-x3-set-numsys},
together with the fundamental region
of the numeration system
on the right.


\FloatBarrier
\subsubsection*{A curve-set of order 36 whose curves have reflection symmetry}

\IMAGE{333333-rotational-manta-t-36-decomp}

Our second curve-set of order 36 was constructed
so that all curves have reflection symmetry.
Its maps are
{\small
\begin{verbatim}
  A |--> A+B-A+B-A+B-A+B-A-C0C0C0C0C+A-C+A-C+A-C+A+B0B0B0B-A0A0A0A0A
  B |--> B+C-B+C-B+C-B+C-B-A0A0A0A0A+B-A+B-A+B-A+B+C0C0C0C-B0B0B0B0B
  C |--> C+A-C+A-C+A-C+A-C-B+C-B+C-B-A+B-A-C+A+B0B0B0B+C-B+C-B+C-B+C+ \
            A-C+A-C+A+B-A+B+C-B-A0A0A-C0C0C0C0C
\end{verbatim}
}
%
The mutual self-similarity of the three curves
is shown in Figure~\ref{fig:333333-rotational-manta-t-36-decomp}.

\IMAGE{333333-rotational-manta-t-36-tile-decomp}
The self-similarity of the CCW tile \tileS{A+B+C+}{1}
is shown in Figure~\ref{fig:333333-rotational-manta-t-36-tile-decomp}.
The tile is a lattice tile and corresponds to the fundamental region
of a numeration system with 36 digits.



\FloatBarrier
\subsubsection*{A curve-set with uneven growth}

Our curve-set has order 6 and maps
\begin{verbatim}
  A |--> A0A+B-A
  B |--> B-A+B+C-B
  C |--> C+A-C0C-B-A+B+C0C
\end{verbatim}
Figure~\ref{fig:333333-rotational-keili} shows
the first, second, and third iterate of the prototile \tileS{A+B+C+}{1}.
So far, nothing suspect.

\IMAGE{333333-rotational-keili}
\IMAGE{333333-rotational-keili-it}

However, at higher iterates an uneven growth shows,
see Figure~\ref{fig:333333-rotational-keili-it}.
The curve-set \emph{is} plane-filling, but not in a nice way:
the growth in one direction is much faster than in the other.
This curve-set was found by Daniel Fischer~\cite{dan-master-thesis}.

%
%
%
%

\subsubsection*{Using the same motif on different colorings}

\IMAGE{33333-motif-col-it1}

Using the same (shape of) motif with different colorings leads
to different shapes of curves (and tiles) in general.
Figure~\ref{fig:33333-motif-col-it1} shows
the CCW tiles \tileS{A+B+C}{1}
with lettering for the coloring in this section (first and second iterate, left)
and lettering according to the coloring of
Section~\jjref{sect:col-333333-star} (again first and second iterate, right).
The coloring used in the images is by curves:
blue for \texttt{A},
red for \texttt{B}, and
green for \texttt{C}.

The seventh iterates of the CCW tile for the coloring in this section
is shown on the left in Figure~\ref{fig:33333-motif-col-itx}
(the CW tile looks identical).
The green part corresponding to the 1-dimensional curve for
\texttt{C} isn't visible at this iterate.
Middle and right respectively show the CCW and CW tiles
for the coloring of
Section~\ref{sect:col-333333-star}.
The green part corresponding to the curve for
\texttt{C} suggests it is 2-dimensional.

\IMAGE{33333-motif-col-itx}

The lozenge-shaped tile may not look terribly exciting
but there is a numeration system
with the fundamental region of just that shape
with radix $-2$ and the set of digits
$\left\{ 0,\, +1,\, +\omega_6,\, +1+\omega_6 \right\}$
where $\omega_6 = \exp(2\pi{}i/6)$.
In addition the shapes of both curves \texttt{A}
and \texttt{B} are triangles, which does not
happen for any curve on the minimal coloring.





\subsubsection{A coloring with four colors}

\IMAGE{333333-bizarro-grid}


The coloring we use is shown in
Figure~\ref{fig:333333-bizarro-grid},
together with its prototiles.
%

Our curve-set of order 9 has maps
\begin{verbatim}
  A |--> A+A-C-D0C+D+B-B0A
  B |--> B+C-D-A0B+C+D-A0B
  C |--> C+D-A-C0D+B+C-D+B-B-B0A+A+A-C
  D |--> D0C0D
\end{verbatim}
The first iterates of the curves are shown in
Figure~\ref{fig:333333-bizarro-curve-r9-var0-it1}.

\IMAGE{333333-bizarro-curve-r9-var0-it1}
\IMAGE{333333-bizarro-curve-r9-var0-it}
Though the motifs of all curves have reflection symmetry,
only curves \texttt{C} and \texttt{D}
keep it with higher iterates.
This can be seen in
Figure~\ref{fig:333333-bizarro-curve-r9-var0-it}
where mutual self-similarity of the curves is also shown.

\IMAGE{333333-bizarro-curve-r9-var0-tiling}
The curves tile the plane as shown in
Figure~\ref{fig:333333-bizarro-curve-r9-var0-tiling}.
Note that the blue regions
each consist of three curves \texttt{A}
and the red regions
each consist of three curves \texttt{B}.


\FloatBarrier

\subsection{The trihexagonal grid}

Here turns are by $60$ degrees ($2\,\pi/6$).

\subsubsection{The minimal coloring with one color}

\IMAGE{3636-grid-graph}
%

Examples of curves on the minimal coloring
of the trihexagonal (or trihex- or $(3.6.3.6)$-grid or kagome lattice,
see Figure~\ref{fig:3636-grid-graph})
are given in \cite[Section~6.3, pp.~66--71]{arndt-curve-search}.

Figure~\ref{fig:3636-unique-r67-1} shows the motif (left)
and both tiles of a curve of order 67 with the map
{\scriptsize
\begin{verbatim}
  F |--> F--F+F--F+F--F+F--F+F+F--F+F--F+F+F--F+F--F+F+F+F--F+F+F--F--F+F+F--F+F--F+F+F+F-- \
           F+F+F--F+F+F+F--F+F+F--F+F--F+F+F+F+F+F--F+F--F+F+F+F--F+F+F+F--F+F+F+F+F
\end{verbatim}
}

\IMAGE{3636-unique-r67-1}

\IMAGE{3636-unique-r67-1-tiles}
We chose this curve because the motif does not look promising,
still, this is a plane-filling curve.
The second iterates of both tiles are shown in
Figure~\ref{fig:3636-unique-r67-1-tiles}.

\FloatBarrier
\subsubsection{A coloring with two colors}

\IMAGE{3636-AB-grid-graph}
The coloring is shown in Figure~\ref{fig:3636-AB-grid-graph},
together with its prototiles
\tileS{A+B+}{6} (CCW, hexagon), \tileS{A--}{3} and \tileS{B--}{3} (CW, triangles).
The CCW tile appears in just one orientation, so it is a lattice tile.

\FloatBarrier
\subsubsection*{A curve-set where the CW prototiles have maximal symmetry}

Our curve-set is of order 25 and has maps
\begin{verbatim}
  A |--> A+B+A--A+B+A--A--A+B+A--A+B--B+A+B+A--A+B+A+B+A+B--B+A--A
  B |--> B+A--A+B+A+B+A+B--B--B+A+B--B--B+A+B--B+A--A+B+A+B+A+B--B
\end{verbatim}

\IMAGE{3636-AB-r25-curve-motifs}

\IMAGE{3636-AB-r25-curve-tile-minus}
\IMAGE{3636-AB-r25-curve-tile-minus-decomp}
The motifs for curves \texttt{A} and \texttt{B} and
CW tiles \tileS{A--}{3} and \tileS{B--}{3}
are shown in Figure~\ref{fig:3636-AB-r25-curve-motifs}.
The third iterates of the CW tiles
are shown in Figure~\ref{fig:3636-AB-r25-curve-tile-minus}.
The symmetry of the CW tiles is maximal, including reflection symmetry
(this is only possible if the order is a square).
The mutual self-similarity of the CW tiles is shown in
Figure~\ref{fig:3636-AB-r25-curve-tile-minus-decomp}.

\IMAGE{3636-AB-r25-curve-tile-plus}
The CCW tile \tileS{A+B+}{3} has 3-fold symmetry,
as shown in Figure~\ref{fig:3636-AB-r25-curve-tile-plus}.

\FloatBarrier
\subsubsection*{A curve-set with one letter constant}

\IMAGE{3636-AB-r16-B-const}

Our curve-set is of order 9 and has the maps
{\scriptsize
\begin{verbatim}
  A |--> A+B+A--A--A+B+A--A+B+A+B+A+B--B+A--A+B--B--B+A+B+A+B--B--B+A+B+A+B+A--A
  B |--> B
\end{verbatim}
}

\IMAGE{3636-AB-r16-B-const-it4}

The first, second, and third iterates of the curve \texttt{A} are shown
in Figure~\ref{fig:3636-AB-r16-B-const}.
Note how the red areas vanish with higher iterates.
The shape of the curve has a 2-fold symmetry in the limit.

\FloatBarrier
The fourth iterate (Figure~\ref{fig:3636-AB-r16-B-const-it4})
is indistinguishable from the curve of order 16
on the triangle grid with map
{\small
\begin{verbatim}
  F |--> F+F-F-F+F-F+F+F0F-F-F+F-F+F+F-F
\end{verbatim}
}
The map for \texttt{F} can be obtained from that of \texttt{A} as follows.
{\scriptsize
\begin{verbatim}
  # map for A:
  A |--> A+B+A--A--A+B+A--A+B+A+B+A+B--B+A--A+B--B--B+A+B+A+B--B--B+A+B+A+B+A--A
\end{verbatim}
}
{\small
\begin{verbatim}
  # drop all (constants) B:
  A |--> A++A--A--A++A--A++A++A+--+A--A+----+A++A+----+A++A++A--A
  # normalize turns
  A |--> A++A--A--A++A--A++A++A0A--A--A++A--A++A++A--A
\end{verbatim}
}
All nonzero turns are by $\pm{}120$,
so we can replace all turns by a single symbol
and use turns by $120$ degrees.
Finally, replace \texttt{A} by \texttt{F}
to arrive at the map given above.

Removing all edges \texttt{B} in the grid
(\jjref{fig:3636-AB-grid-graph})
gives the triangle grid with minimal coloring, so this is
not too surprising.

\FloatBarrier
\subsubsection{A coloring with three colors}


\IMAGE{3636-FGH-grid-graph}
We use the coloring shown in Figure~\ref{fig:3636-FGH-grid-graph}.
The CCW prototile \tileS{F+G+H+}{2} now has 2-fold symmetry
and the CW prototile \tileS{F--G--H--}{1} appears in two orientations.

With this coloring the sequence of letters in any curve is always
\texttt{F}, \texttt{G}, \texttt{H}, \texttt{F}, \texttt{G}, \texttt{H}, \ldots,
regardless of the turns.

\IMAGE{3636-FGH-motifs}
We use a curve-set of order 25 with maps
{\small
\begin{verbatim}
  F |--> F+G+H--F+G+H--F--G+H--F+G+H+F+G--H--F+G+H+F--G+H+F+G+H+F--G+H--F
  G |--> G+H+F+G--H--F+G+H+F--G--H+F+G--H+F--G+H+F+G--H+F+G+H+F+G--H+F--G
  H |--> H+F+G--H+F+G--H--F+G+H--F+G--H+F+G+H+F+G--H
\end{verbatim}
}
The motifs of the curves are shown in
Figure~\ref{fig:3636-FGH-motifs}.

\IMAGE{3636-FGH-tile-minus}
The self-similarity of the CW prototile
is shown in Figure~\ref{fig:3636-FGH-tile-minus},
note that the red and blue small copies appear in two orientations.

\IMAGE{3636-FGH-tile-plus}
Figure~\ref{fig:3636-FGH-tile-plus} shows the 2-symmetric CCW prototile.
It is a lattice tile.

\FloatBarrier
\subsubsection{A coloring with six colors}

%
In the last few sections we have reduced the rotational symmetry of
the CCW tile from 6-fold, to 3-fold and 2-fold.
A coloring with six colors is shown in
Figure~\ref{fig:3636-set-grid},
it guarantees no rotational symmetry for the CCW tile.
It still is a lattice tile.

\IMAGE{3636-set-grid}

\IMAGE{3636-set-curve-r16-var0}
Our curve-set is of order 16 and has maps
{\small
\begin{verbatim}
  A |--> A--E+F+A+B+C+D--B+C--A+B+C--A--E+F--D+E+F+A+B+C--A--E+F+A
  B |--> B+C+D--B--F+A+B
  C |--> C--A+B+C+D+E+F--D+E--C+D--B--F+A+B+C+D+E--C--A+B+C
  D |--> D--B+C+D+E+F+A--E+F--D--B+C+D
  E |--> E--C+D+E+F+A+B--F+A--E+F--D--B+C+D+E+F+A--E--C+D+E
  F |--> F+A+B--F--D+E+F
\end{verbatim}
}
The motifs of the curves, assembled into the CCW tile,
are shown in
Figure~\ref{fig:3636-set-curve-r16-var0}.

\IMAGE{3636-set-curve-r16-CW-tiles}
The CW tiles are shown in
Figure~\ref{fig:3636-set-curve-r16-CW-tiles}.
There are two ways of pairing them into lattice tiles.
The first way is to connect them so that
the turn \texttt{B+C} is contained.
This arrangement
and its self-similarity is shown in
Figure~\ref{fig:3636-set-curve-r16-lattice-tile1}.
\IMAGE{3636-set-curve-r16-lattice-tile1}
\IMAGE{3636-set-curve-r16-lattice-tile2}

The second way to connect the tiles contains
the turn \texttt{A+B}.
This arrangement
is shown in
Figure~\ref{fig:3636-set-curve-r16-lattice-tile2}.
It is not self-similar:
note that the small copy
at the  top of the right image in the figure
is incomplete, it is missing the red/yellow/pink part
which lies at the lower right of the image.


\FloatBarrier
\clearpage

\subsection{The \texorpdfstring{$(3.4.6.4)$}{(3.4.6.4)}-grid}

Here we use turns by $30$~degrees ($2\,\pi/12$).

\subsubsection{The minimal coloring with two colors}

\IMAGE{3464-tile-repeat}
The minimal coloring of the $(3.4.6.4)$-grid
and its prototiles are shown in
Figure~\ref{fig:3464-tile-repeat}
(same as Figure~\jjref{fig:3464-tile}).
%

A tedious construction for edge-covering curves on this grid
is given in \cite[Section~4.2.1, pp.42-43]{arndt-curve-search}.
Using the proper amount of edge classes (two, not one!)
makes life considerably easier.

\FloatBarrier
\subsubsection*{6-fold symmetry is possible for the CCW tile \texorpdfstring{\tileS{B++++}{3}}{[B+++]**3}}

\IMAGE{3464-curve-r9-var1-motifs}
The rotational symmetries of the prototiles are
six-fold for the \tileS{A++}{6},
at least three-fold for \tileS{B++++}{3},
and two-fold for \tileS{A---B---}{2}.
The prototile \tileS{B++++}{3}  may have six-fold symmetry.
An example is the curve-set of order 9 with maps
{\small
\begin{verbatim}
  A |--> A++A++A++A++A---B---A---B++++B---A
  B |--> B++++B---A---B++++B++++B---A---B
\end{verbatim}
}
The motifs of the curves and tiles are shown in
Figure~\ref{fig:3464-curve-r9-var1-motifs}.

\IMAGE{3464-curve-r9-var1-koch-tiles}
\IMAGE{3464-curve-r9-var1-koch-tiles-decomp}
The first iterate of \tileS{B++++}{3} has six-fold symmetry,
see Figure~\ref{fig:3464-curve-r9-var1-koch-tiles}.
Both prototiles happen to have the shape of the Koch snowflake,
see \cite[Figure~3]{burns-fractal-tilings}.
Note that this shape also has reflection symmetry.

\IMAGE{3464-curve-r9-var1-tile-minus}
The tile \tileS{A---B---}{2}
shown in Figure~\ref{fig:3464-curve-r9-var1-tile-minus}
can be decomposed into 9 smaller copies of itself,
appearing in three orientations (image on the right).

\FloatBarrier

\IMAGE{3464-curve-r39-hex-var1-tiles-plus}
The two CCW tiles shown in
Figure~\ref{fig:3464-curve-r39-hex-var1-tiles-plus}
also have 6-fold symmetry,
they belong to a curve-set of order 39.
Here the shapes of the two tiles are different from each other.
The maps of this curve-set, omitting the turns, are
{\small
\begin{verbatim}
  A |--> AABAAAAAABBABBABBABAAAABBBABAABBAAAAAABA
  B |--> BBAAABABBABBBABBBAB,BBAAABABBABBBABBBAB
\end{verbatim}
}
Note that the map for \texttt{B} is a square,
as indicated by the comma in its middle.
Here we can omit turns, as the transitions are
\texttt{A++A}, \texttt{B++++B}, \texttt{A---B}, and \texttt{B---A}
and the two letters of each transition uniquely determine the turn
(this is of course \emph{not} true in general).

\FloatBarrier
\subsubsection*{A curve-set of order 37}

The following curve-set is presented because
the tiles look quite curious.
It has order 37 and the maps
(omitting turns again)
{\small
\begin{verbatim}
  A |--> AAAABAAABBABAAABBAA
  B |--> BABBBABAAABBBABBBABBBABAABBBABBBABBBAABABBBAAABAAABAAAB
\end{verbatim}
}

\IMAGE{3464-r37-tiles-plus}
\IMAGE{3464-r37-tile-plus-A6-decomp}
\IMAGE{3464-r37-tile-plus-B3-decomp}
The CCW tiles \tileS{A++}{6} and \tileS{B++++}{3} are
shown in Figure~\ref{fig:3464-r37-tiles-plus}.
Their mutual self-similarity can be observed
in Figures~\ref{fig:3464-r37-tile-plus-A6-decomp} and \ref{fig:3464-r37-tile-plus-B3-decomp}.
Note that while the prototiles \tileS{A++}{6} are separated,
the prototiles \tileS{B++++}{3} are not,
making the images slightly difficult to interpret.
The gray borders drawn between the curves are to mitigate this.


\FloatBarrier
\subsubsection*{A curve-set with letter \texttt{A} constant}

Curve-sets whose L-system leaves some of the letters constant exist.
Here we have just one non-constant letter and
the curve for it shows simple self-similarity.

\IMAGE{3464-degenerate-curve-A-const-motif}
\IMAGE{3464-degenerate-curve-A-const-decomp}

\FloatBarrier
Figure~\ref{fig:3464-degenerate-curve-A-const-motif} shows the
motif and second iterate of curve \texttt{B}
of a curve-set of order 9 with maps
{\small
\begin{verbatim}
  A |--> A
  B |--> BAAAAABABBBABBABB
\end{verbatim}
}
Figure~\ref{fig:3464-degenerate-curve-A-const-decomp}
shows the self-similarity of the curve.
Note that the curve has the shape of half of
the curve shown on the right in
Figure~\jjref{fig:3464-curve-r9-var1-tile-minus}.

\subsubsection*{A curve-set with letter \texttt{B} constant}

\IMAGE{3464-degenerate-curve-B-const-motif}

Figure~\ref{fig:3464-degenerate-curve-B-const-motif} shows the
motif and second iterate of curve \texttt{A}
of a curve-set of order 19 with maps
{\small
\begin{verbatim}
  A |--> AABBAAABABBABBBAABABBBABBBABAAAAAABBA
  B |--> B
\end{verbatim}
}

\IMAGE{3464-degenerate-curve-B-const-decomp}
%
The self-similarity of the curve is shown in
Figure~\ref{fig:3464-degenerate-curve-B-const-decomp}.
The apparent lozenges (for example, at the top right)
are the tiles \tileS{A---B---}{2}, colored
by their orientation.
The tile \tileS{A++}{6} (6-fold symmetry) can be seen
at the lower left.

\FloatBarrier
\subsubsection{A coloring with four colors}

\IMAGE{3464-rotational-grid-graph}

The coloring shown in Figure~\ref{fig:3464-rotational-grid-graph}
splits the triangular CCW prototile of the minimal coloring
into two classes, prototiles
\tileS{B++++}{3} and \tileS{b++++}{3}.

\subsubsection*{First curve-set}

The curve-set of order 9 has maps
{\small
\begin{verbatim}
  A |--> A---B++++B++++B---a++A---B---a---b++++b---A++a++A
  a |--> a++A++a---b---A---B++++B---a---b++++b++++b---A++a
  B |--> B++++B---a++A---B
  b |--> b++++b---A++a---b
\end{verbatim}
}
\IMAGE{3464-rotational-x00-motifs}

\IMAGE{3464-rotational-x00-prototiles-plus}
\IMAGE{3464-rotational-x00-prototiles-plus-monochrome}
The CCW tiles shown in
Figures~\ref{fig:3464-rotational-x00-prototiles-plus}
and~\ref{fig:3464-rotational-x00-prototiles-plus-monochrome}
have reflection symmetry.
Moreover the tile \tileS{A++a++}{3}
has 6-fold rotational symmetry.

\FloatBarrier
\subsubsection*{Second curve-set}

\IMAGE{3464-rotational-x02-motifs}
The curve-set of order 9 has maps
{\small
\begin{verbatim}
  A |--> A---B++++B++++B---a++A---B---a---b++++b---A++a++A
  a |--> a++A++a---b---A---B++++B---a++A++a
  B |--> B++++B---a---b++++b++++b---A---B
  b |--> b++++b---A++a---b
\end{verbatim}
}

\IMAGE{3464-rotational-surgery}
This curve-set was created by modifying the interior of
the CW prototile as shown in Figure~\ref{fig:3464-rotational-surgery}.
This technique of removing a polygon from one curve and reattaching it to
another is quite useful for modifying curve-sets.

\IMAGE{3464-rotational-x02-prototiles-plus}
\IMAGE{3464-rotational-x02-prototiles-plus-monochrome}

\IMAGE{3464-rotational-x02-prototiles-minus}

The CCW tiles again have reflection symmetry, they are shown in
Figures~\ref{fig:3464-rotational-x02-prototiles-plus}
and~\ref{fig:3464-rotational-x02-prototiles-plus-monochrome}.

The CW prototile is shown in
Figure~\ref{fig:3464-rotational-x02-prototiles-minus}.
It has the same shape as the CW prototile for the
previous section (we did the change inside it).

\FloatBarrier
\subsubsection{A coloring with six colors}

\IMAGE{3464-fhg-grid-graph}
The coloring shown in Figure~\ref{fig:3464-fhg-grid-graph}
has three CW prototiles,
\tileS{F---f---}{2},
\tileS{G---g---}{2}, and
\tileS{H---h---}{2}.
Each of these has 2-fold rotational symmetry
and appears in just one orientation.
%

%


Our curve-set of order 9 has maps
{\small
\begin{verbatim}
  F |--> F++G++H++F++G---g++++f---F---f---F++G---g++++f---F
  f |--> f++++h---H---h++++g++++f---F---f
  G |--> G++H++F++G++H---h---H---h++++g---G
  g |--> g++++f---F---f++++h++++g---G---g
  H |--> H++F++G++H---h---H
  h |--> h++++g---G---g++++f++++h---H---h
\end{verbatim}
}

\IMAGE{3464-fhg-r9-var0-tile-minus-decomp}
The mutual self-similarity of the three tiles is shown in
Figure~\ref{fig:3464-fhg-r9-var0-tile-minus-decomp}.
Observe how each tile appears in just one orientation.
Borders are drawn between the curves, this helps to
spot the borders between the small copies of
the blue prototile on the left.
%


\FloatBarrier
\subsection{The \texorpdfstring{$(3^6;\, 3.3.6.6)$}{(3.3.3.3.3.3; 3.3.6.6)}-grid}

Here we use turns by $60$ degrees ($2\,\pi/6$).

\IMAGE{333333-3366-grid-graph}
The $(3^6; 3.3.6.6)$-grid
is shown in Figure~\ref{fig:333333-3366-grid-graph}.
The grid has five classes of edges,
denoted by the letters \texttt{A} \ldots \texttt{E}.
The transitions are\begin{verbatim}
  A+C,  B++D, C+A,   D++E,  E++B,
  A-B,  B-A,  C--D,  D--E,  E--C,
  D0E
\end{verbatim}
Note that there are three transitions from \texttt{D} to \texttt{E}.
The CCW prototiles are \tileS{A+C+}{3} and \tileS{B++D++E++}{1},
and the CW prototiles are \tileS{A-B-}{3} and \tileS{C--D--E--}{1}.
We only use this minimal coloring of the grid.

\subsubsection*{A curve set of order 16}

Our curve-set has order 16 and maps
{\small
\begin{verbatim}
  A |--> A-B-A+C+A+C+A+C--D++E--C+A-B-A-B-A-B++D++E--C--D++E++B-A
  B |--> B-A-B++D--E++B++D--E++B++D--E--C+A-B
  C |--> C+A+C--D++E--C--D++E--C--D++E++B-A+C
  D |--> D++E--C--D0E++B++D--E++B-A-B-A-B++D
  E |--> E--C+A+C+A+C--D++E--C--D0E++B++D--E
\end{verbatim}
}
It was designed so that the curves have mostly straight borders.

The mutual self-similarity of the five curves
is shown in Figure~\ref{fig:333333-3366-r16-var0-curve-decomp}.
The shapes of the CCW prototiles are shown
in Figure~\ref{fig:333333-3366-r16-var0-tiles-plus}.
The tiling of the plane by the five curves is shown
in Figure~\ref{fig:333333-3366-r16-var0-tiling}.

\IMAGE{333333-3366-r16-var0-curve-decomp}

\IMAGE{333333-3366-r16-var0-tiles-plus}

\IMAGE{333333-3366-r16-var0-tiling}


\FloatBarrier
\subsubsection*{A curve set of order 13 with just one non-constant symbol}


\IMAGE{333333-3366-trans-r13-15}
The curve-set of order 13 has the following map for the letter \texttt{A},
all other letters are constants.
{\scriptsize
\begin{verbatim}
  A |--> A-B++D++E--C+A-B++D0E--C+A-B++D0E--C+A-B++D--E--C+A-B++D--E--C+A-B++D++E--C+A-B++D0E \
         --C+A-B++D++E--C+A-B++D++E--C+A-B++D--E--C+A-B++D0E--C+A-B++D--E--C+A
\end{verbatim}
}

\IMAGE{333333-3366-trans-r13-15-curve-decomp}

It was created to look similar to the curve of order 13 named \texttt{R13-15}
from \cite[Figure~3.1-A, p.~15]{arndt-curve-search}
with map
\begin{verbatim}
  F |--> F+F0F0F-F-F+F0F+F+F-F0F-F
\end{verbatim}
For higher iterates the two curves are hard to tell apart
as their limiting shapes are identical.
The fourth iterate is shown
in Figure~\ref{fig:333333-3366-trans-r13-15-curve-decomp}.
%
%
Iterates of the tile \tileS{A+C+}{3}
are shown in Figure~\ref{fig:333333-3366-trans-r13-15-tiles-plus}.
The tile has 6-fold rotational symmetry in the limit.
For higher iterates all but the blue and green regions
(for the letter \texttt{A}) will vanish.

\IMAGE{333333-3366-trans-r13-15-tiles-plus}

Indeed any curve on the minimal coloring of the triangle grid can be turned to
a curve-set on the $(3^6;\, 3.3.6.6)$-grid
with all letters but \texttt{A} constant:
In the map for the curve on the triangle grid,
replace
all \texttt{+} by \texttt{p},
all \texttt{-} by \texttt{m},
all \texttt{0} by \texttt{n},
all \texttt{p} by \texttt{B++D++E--C},
all \texttt{m} by \texttt{B++D--E--C},
all \texttt{n} by \texttt{B++D0E--C}, and
all \texttt{F} by \texttt{+A-},
finally remove the leading \texttt{+} and the trailing \texttt{-}.
This gives the map for \texttt{A}.
In our example, the following command produces
the map for \texttt{A} shown above (split into pieces for readability).
{\small
\begin{verbatim}
echo 'F+F0F0F-F-F+F0F+F+F-F0F-F' | \
  sed 's/+/p/g; s/-/m/g; s/0/n/g;' | \
  sed 's/p/B++D++E--C/g; s/m/B++D--E--C/g; s/n/B++D0E--C/g; s/F/+A-/g;' | \
  sed 's/^+//; s/-$//'
\end{verbatim}
}

Note that dropping all but the \texttt{A}-edges in the grid
(Figure~\jjref{fig:333333-3366-grid-graph})
gives the triangle grid with minimal coloring.

%
%
%


\FloatBarrier

\subsection{The \texorpdfstring{$(3.4.4.6;\, 3.4.6.4)$}{(3.4.4.6; 3.4.6.4)}-grid}


\IMAGE{3446-3464-grid-graph-repeat}
The $(3.4.4.6; 3.4.6.4)$-grid and its prototiles are
shown in Figure~\ref{fig:3446-3464-grid-graph-repeat}
(same as Figure~\jjref{fig:3446-3464-grid-graph}).
The prototiles are
\tileS{A++}{6},
\tileS{B++++C++++D++++}{1},
\tileS{E+++F+++}{2},
\tileS{A---B---F---D---}{1},
and
\tileS{C--E--}{3}.

Turns are by 30 degrees $(2\,\pi/12)$.
%
The transitions are
\texttt{A++A},
\texttt{A---B},
\texttt{B++++C},
\texttt{B---F},
\texttt{C++++D},
\texttt{C--E},
\texttt{D++++B},
\texttt{D---A},
\texttt{E+++F},
\texttt{E--C}, and
\texttt{F+++E}, and
\texttt{F---D}.
As the pair of letters at start and end uniquely determine
the turn, we can omit the turns in the maps to save space.

\FloatBarrier
\subsubsection*{A curve-set of order 7 with four constant letters}
We give an example of a plane-filling curve-set where
only two of the five prototiles are plane-filling.
It has order 7 and just two non-constant maps:
{
\begin{verbatim}
  A |--> AAABFEFDABCECECDBFEFDAAA
  F |--> FDBFDBFECDBFECECECDBF
\end{verbatim}
}
The maps for the letters
\texttt{B}, \texttt{C}, \texttt{D}, and \texttt{E} are constant.
The motifs for \texttt{A} and \texttt{F} are shown in
Figure~\ref{fig:3446-3464-curve-set-var-2-motifs}.
\IMAGE{3446-3464-curve-set-var-2-motifs}

%
%
\IMAGE{3446-3464-curve-set-var-2-prototile-A6}
\IMAGE{3446-3464-curve-set-var-2-prototile-EF2}
\IMAGE{3446-3464-curve-set-var-2-prototile-ABFD1}
\FloatBarrier
Iterates of the prototiles
\tileS{A++}{6}, \tileS{E+++F+++}{2}, and \tileS{A---B---F---D---}{1}
are respectively shown in
Figures
\ref{fig:3446-3464-curve-set-var-2-prototile-A6},
\ref{fig:3446-3464-curve-set-var-2-prototile-EF2}, and
\ref{fig:3446-3464-curve-set-var-2-prototile-ABFD1}.
The prototiles \tileS{B++++C++++D++++}{1} and \tileS{C--E--}{3}
have all edges constant.

\IMAGE{3446-3464-decomp-var-2}
The mutual self-similarities of the curves \texttt{A} and \texttt{F}
are shown in Figure~\ref{fig:3446-3464-decomp-var-2}.
The substitution matrix for this curve-set is
\begin{eqnarray*}
M  & = &
\begin{bmatrix}
7 & 0 & 0 & 0 & 0 & 0 \\
2 & 1 & 0 & 0 & 0 & 4 \\
2 & 0 & 1 & 0 & 0 & 4 \\
2 & 0 & 0 & 1 & 0 & 4 \\
2 & 0 & 0 & 0 & 1 & 4 \\
2 & 0 & 0 & 0 & 0 & 5 \\
\end{bmatrix}
\end{eqnarray*}
The curve for \texttt{A} (corresponding to the first column of the matrix)
has dimension $2$.
The curve for \texttt{F} (last column)
has dimension $2\cdot\log_7{5}\approx{}1.65417$,
as the map for \texttt{F} contains five letters \texttt{F}
(the lower right entry in the matrix)
and no letters corresponding to curves with dimension 2.
The constant letters \texttt{B}, \texttt{C}, \texttt{D}, and \texttt{E}
do not contribute to the dimensions
of \texttt{A} and \texttt{F}.
Any map containing letters corresponding to a curve with dimension 2
gives a 2-dimensional curve.
As the curve for \texttt{A} is the only 2-dimensional curve in this curve-set,
the map for \texttt{F} cannot contain the letter \texttt{A}.
It indeed does not, see the upper right entry of the matrix.

\subsubsection*{A curve-set of order 31}

\IMAGE{3446-3464-r31-var0-motifs}
%
Our curve-set is of order 31 and has maps
\begin{verbatim}
  A |--> ABCDAABFEFDABCECECDBFEFDABCECDBFDBFDBCEFDBCECEFDABCDA
  B |--> BFEFDBFECECECDBFDAAAABCDAAB
  C |--> CEFDABCEFEFDAABCDAAABFEC
  D |--> DBFDAAAAAABFDBCECECEFEFD
  E |--> ECDABFECDBFDBFEFDBCECECDABCDAABFDBFE
  F |--> FECECECDBFEFDABCDABFEF
\end{verbatim}
The first iterates of the curves are shown in
Figure~\ref{fig:3446-3464-r31-var0-motifs}.

\IMAGE{3446-3464-grid-mod-180deg}
A coloring by orientation of the \emph{undirected} edges
of the grid is shown in
Figure~\ref{fig:3446-3464-grid-mod-180deg}.
This coloring brings out the bands of squares in the grid.
For the CW tiles in
Figure~\ref{fig:3446-3464-r31-var0-tiles-minus-2}
this coloring is used.
%
\IMAGE{3446-3464-r31-var0-tiles-minus-2}

\IMAGE{3446-3464-r31-var0-tiles-plus}
The mutual self-similarities of
the CCW tiles is shown in
Figure~\ref{fig:3446-3464-r31-var0-tiles-plus}.
Line drawings of the motifs on the left,
area drawings of the third iterates on the right.

\FloatBarrier
\clearpage

\section{Curve-sets for grids with odd valencies}\label{sect:double-edges}


\IMAGE{hexagon-grid}
So far we gave curve-sets for grids where every
point has even valency (number of incident edges).

In the hexagon grid shown in
Figure~\ref{fig:hexagon-grid}
every point has valency 3.
Every time a curve visits a point it needs to leave.
That is, we seemingly have the problem that at a point
of odd valency we cannot avoid creating a dead end.
%

The solution is quite simple: edge bifurcation.
We declare some edges to be
anti-parallel pairs of edges.
Once all points have even valency we
can work (almost) as before.
For our examples we always turn all edges
into double edges.

Any such pair of \jjxterm{double edges}{double edge}
creates a \jjterm{digon},
which has to be considered a prototile.
At both ends of a digon a turn by 180 degrees (U-turn)
can appear, we reserve the symbol \texttt{!} to specify
U-turns in the Lindenmayer system.
In our examples these digon tiles are always CW.

One can also create double edges where it is not strictly needed,
we will do so for the square grid and the triangle grid.

\IMAGE{d333333-curve-r13-1-solid}

In the area drawings the area assigned to an edge will be on its
left, see Figure~\ref{fig:d333333-curve-r13-1-solid}
(compare to Figure~\jjref{fig:333333-curve-r13-1-solid}).


\FloatBarrier
\clearpage
\subsection{The hexagon grid}

\IMAGE{d666-grid-graph}
In the hexagon grid with double edges
is shown in Figure~\ref{fig:d666-grid-graph}.
Each edge is responsible for an area to its left,


There is just one edge class;
we use the letter \texttt{A} for it.

Here we use turns by 60 degrees ($2\pi/6$).

The transitions are
\texttt{A+A}, \texttt{A-A}, and \texttt{A!A}.
There is one CCW prototile \tileS{A+}{6},
a lattice tile,
and one CW prototile \tileS{A!}{2},
a digon tile.

%
\subsubsection*{A curve of order 25}

\IMAGE{d666-r25-var1-outline}

Creating a curve is quite easy.
Start with a (self-avoiding) curve with 2-fold symmetry.
The symmetry is necessary to make the digon tile work.

Now we have a CCW tile \tileS{A+}{6} which is not yet filled,
shown in Figure~\ref{fig:d666-r25-var1-outline} (left).
The unfinished curve has the (incomplete) map
\begin{verbatim}
  A |--> A+A-A+A+A+A+A-A-A-A-A+A-A
\end{verbatim}
Filling the tile is possible in more than one way.
The completion shown on the right corresponds
to a curve with the map
\begin{verbatim}
  A |--> A+A+A!A+A+A+A+A+A-A+A!A+A+A!A+A+A+A+A!A-A-A+A+A-A
\end{verbatim}

\IMAGE{d666-r25-var1-tile-solid}
\IMAGE{d666-r25-var1-digon-solid}
Area drawings of the tiles are shown in
Figure~\ref{fig:d666-r25-var1-tile-solid} (CCW)
and
Figure~\ref{fig:d666-r25-var1-digon-solid} (CW, digon).


\FloatBarrier
\subsubsection*{A curve of order 109}

\IMAGE{d666-r109-var1-motif}

The motif of a curve of order 109 with map
{\scriptsize
\begin{verbatim}
  A |--> A+A+A!A+A+A+A!A+A+A+A!A+A+A+A!A+A+A+A+A+A!A+A+A+A!A+A+A+A!A+A+A+A!A+A+A+ \
         A+A!A+A+A+A!A+A+A+A+A+A!A+A+A+A!A+A+A+A+A+A-A+A!A+A-A-A+A-A+A-A+A!A+A-A+ \
         A-A+A-A+A!A+A+A+A!A+A+A+A!A+A+A+A!A+A+A+A!A+A+A+A!A-A+A-A+A-A+A-A+A-A+A-A
\end{verbatim}
}
is shown in Figure~\ref{fig:d666-r109-var1-motif}


\IMAGE{d666-r109-var1-tile-plus-it1}
\IMAGE{d666-r109-var1-tile-plus}

The CCW tile is shown in Figures
 \ref{fig:d666-r109-var1-tile-plus-it1}
and
 \ref{fig:d666-r109-var1-tile-plus},
the CW (digon) tile in
Figure~\ref{fig:d666-r109-var1-tile-digon}.

\IMAGE{d666-r109-var1-tile-digon}

\FloatBarrier
\subsubsection*{A curve of order 19}

\IMAGE{d666-r19-it1}

Figure~\ref{fig:d666-r19-it1}
shows the motif and second iterate of a curve of order 19
with map
\begin{verbatim}
   A |--> A+A-A+A+A!A-A+A+A!A+A+A!A+A+A+A+A+A-A
\end{verbatim}

\IMAGE{d666-r19-it}
The fourth iterate is shown in Figure~\ref{fig:d666-r19-it}.
Now compare to
Figure~\jjref{fig:3464-degenerate-curve-B-const-decomp},
the images look identical.

Indeed, take the map for the non-constant letter \texttt{A}
for the curve of Figure~\ref{fig:3464-degenerate-curve-B-const-decomp},
{\small
\begin{verbatim}
   A |--> A++A---B++++B---A++A++A---B---A---B++++B---A---B++++B++++B---A++ \
          A---B---A---B++++B++++B---A---B++++B++++B---A---B---A++A++A++A++ \
          A++A---B++++B---A
\end{verbatim}
}
Remove all letters \texttt{B}:
{\small
\begin{verbatim}
  A |--> A++A---++++---A++A++A------A---++++---A---++++++++---A++A------ \
         A---++++++++---A---++++++++---A------A++A++A++A++A++A---++++---A
\end{verbatim}
}
Normalize turns:
{\small
\begin{verbatim}
  A |--> A++A--A++A++A------A--A++A++A------A++A++A------A++A++A++A++A++A--A
\end{verbatim}
}
These turns are by 30 degrees, we want to use turns by 60 degrees:
{\small
\begin{verbatim}
  A |--> A+A-A+A+A---A-A+A+A---A+A+A---A+A+A+A+A+A-A
\end{verbatim}
}
Finally replace \texttt{---} by \texttt{!}, they are U-turns.
This gives the map of our curve.

This is not a surprise.
Reducing the edges \texttt{B} in the $(3.4.6.4)$-grid
(the red ones in Figure~\jjref{fig:3464-tile-repeat})
to zero gives the hexagon grid with double edges.

The lozenge shaped  areas in
Figure~\ref{fig:d666-r19-it}
are the digon tiles \tileS{A!}{2}.

\subsubsection*{A coloring with three edge classes}

\IMAGE{d666-ABC-grid}
We use the coloring with three colors (edge classes) shown in
Figure~\ref{fig:d666-ABC-grid}.
The CCW tile is \tileS{A+B+C+}{2}, it is a lattice tile.
The CW tiles are \tileS{A!}{2}, \tileS{B!}{2}, and \tileS{C!}{2}.

\IMAGE{d666-ABC-r13-var01-motifs}
%
The motifs for a curve-set of order 13
with maps
\begin{verbatim}
  A |--> A+B+C+A+B!B+C+A-C!C+A-C+A+B-A
  B |--> B+C+A+B+C!C+A+B!B-A+B+C-B
  C |--> C+A+B+C+A!A+B!B+C+A-C
\end{verbatim}
is shown in Figure~\ref{fig:d666-ABC-r13-var01-motifs}.

\IMAGE{d666-ABC-r13-var01-tiles-minus}
The mutual self-similarity of the CW tiles is shown in
Figure~\ref{fig:d666-ABC-r13-var01-tiles-minus}.


\FloatBarrier
\subsection{The square grid}\label{sect:d4444}

\IMAGE{d4444-grid-graph}
Though not needed, we can use double edges for the square grid,
see Figure~\ref{fig:d4444-grid-graph}.

There is just one edge class.

Here we use turns by 90 degrees ($2\pi/4$).

The transitions are
\texttt{A+A},
\texttt{A0A},
\texttt{A-A}, and
\texttt{A!A}.
The CCW prototile is \tileS{A+}{4}, it is a lattice tile.
The CW (digon) tile is \tileS{A!}{2}.


\subsubsection*{A curve of order 4}

\IMAGE{d4444-r4-tile-it1}
The only curve of order 4 is well-known,
see
\cite[Section~4.5, pp.60-61]{sagan}
or
\cite[Item~5.3, Section~A.5, p.243]{bader-sfc}.
%
We start with an empty $2\times{}2$ square,
see Figure~\ref{fig:d4444-r4-tile-it1} (left).
There is just one way to fill the square (middle).
The map of the resulting curve is
\begin{verbatim}
  A |--> A+A!A+A
\end{verbatim}
The second iterate of this tile \tileS{A+}{4}
is shown on the right.
The limiting shape of the curve is a triangle.

\IMAGE{d4444-r4-pc-curve}
The curve can be rendered so that it traverses the points
of the $(4.8.8)$-grid as shown in
Figure~\ref{fig:d4444-r4-pc-curve}
and \cite[Figure~5, p.~428]{sagan-reflections}.
%
Here is a recipe for this rendering from our curve:
in the word produced by the Lindenmayer system replace,
in this order,
all \texttt{+} with \texttt{+F+},
all \texttt{A} with \texttt{+F-F-F+},
finally
all \texttt{!} with \texttt{--F--}.
For drawing the curve, use turns by 45 degrees ($2\pi/8$).


\subsubsection*{A curve of order 5}

\IMAGE{d4444-r5-tile-it1}
A curve of order 5 with map
\begin{verbatim}
  A |--> A+A0A!A+A
\end{verbatim}
was given in
\cite[lower left in Figure~2.5, p.~16]{prusinkiewicz-lsys-fract-plants}.
The first iterate of its CCW tile \tileS{A+}{4}
is shown in Figure~\ref{fig:d4444-r5-tile-it1}.

\IMAGE{d4444-r5-tiles}
\IMAGE{r5-dragon}
The sixth iterates of both tiles
\tileS{A+}{4} and \tileS{A!}{2}
are shown in Figure~\ref{fig:d4444-r5-tiles}.
Note that the shape of the tile \tileS{A!}{2}
is that of the curve of order 5
with map
\begin{verbatim}
  F |--> F+F+F-F-F
\end{verbatim}
on the square grid without double edges, see
Figure~\ref{fig:r5-dragon}.
This curve is the only one of order 5
on the square grid with minimal coloring.

\FloatBarrier
\subsubsection*{A curve of order 17}

\IMAGE{d4444-r17-tile-it1}
Creating curves manually is as straightforward
as for the hexagon grid, see
Figure~\ref{fig:d4444-r17-tile-it1}
for a curve of order 17 with map
\begin{verbatim}
  A |--> A+A!A+A+A0A0A!A0A+A+A+A!A+A!A0A+A
\end{verbatim}

\IMAGE{d4444-r17-tiles}
Fourth iterates of the tiles and the self-similarity of the tile \tileS{A!}{2}
are shown in Figure~\ref{fig:d4444-r17-tiles}.
That the small copies on the right appear in two orientations should be no
surprise, \tileS{A!}{2} lies along the edge pairs and those appear in horizontal
and vertical orientation.

\FloatBarrier
\subsubsection*{A curve of order 85}

\IMAGE{d4444-r85-var1-motif}

This curve of order 85
has the map
{\scriptsize
\begin{verbatim}
  A |--> AAAAAA+A!AA+A!A+A+A!A+A+A!A+A+A!A+A+A!A+A!A+A+A!A+A+A!A+A+A!A+A+A!A+A!A+A+ \
         A!A+A+A!A+A+A!A+A!A+A+A!A+A+A!A+A!A+A+A!A+A!A+A+A+AA+AAA+AAAA+AAAAA+AAAAAA
\end{verbatim}
}
The motif is shown in
Figure~\ref{fig:d4444-r85-var1-motif}.

The self-similarity of the tile \tileS{A+}{4} is shown in
Figure~\ref{fig:d4444-r85-var1-tile-plus}.

\IMAGE{d4444-r85-var1-tile-plus}


\FloatBarrier
\subsection{The triangle grid}

%
The triangle grid with double edges is shown in
Figure~\ref{fig:d333333-grid-graph}.
There are now six orientations of the edges in the grid.

Turns are by 60 degrees ($2\pi/6$)
and all six transitions are possible,
\texttt{A0A},
\texttt{A+A},
\texttt{A++A},
\texttt{A!A},
\texttt{A-A}, and
\texttt{A--A}.

There is just one edge class.
The CCW tile is \tileS{A++}{3},
the CW tile is the digon \tileS{A!}{2}.

\IMAGE{d333333-grid-graph}

\FloatBarrier
\subsubsection*{A curve of order 4}

\IMAGE{d333333-szilard-quinton}
Our first example appears (in some disguise) in
\cite[Figure~2.9, p.~19]{prusinkiewicz-lsys-fract-plants},
where it is attributed to \cite{szilard-quinton}.
The maps of the curve given are
\begin{verbatim}
  A |--> AF+F+AF-F-F-AF-F+F+F-F+F+F-A
  F |--> F
\end{verbatim}
The letters \texttt{A} are removed before drawing.
We remove all letters \texttt{F} and normalize the turns.
to obtain the arguably most simple curve on this grid,
see Figure~\ref{fig:d333333-szilard-quinton}.
\begin{verbatim}
  A |--> A++A---A-++-++-A  (drop letters F)
  A |--> A++A!A+A          (normalize turns)
\end{verbatim}

\IMAGE{d333333-szilard-quinton-decor}
Renderings as suggested by Szilard and Quinton are shown in
Figure~\ref{fig:d333333-szilard-quinton-decor}.
The curves cover all points on the hexagon grid.
The axiom for the tile on the right is
\tileS{AF+F+}{3}.

\subsubsection*{A curve of order 3}

\IMAGE{d333333-terdragon-1}
This curve has the map
\begin{verbatim}
  A |--> A++A--A
\end{verbatim}
which is just that of the terdragon,
whose map is usually given as
\begin{verbatim}
  F |--> F+F-F
\end{verbatim}
for turns by 120 degrees ($2\pi/3$).

Here we took \texttt{A+A-A} which has 2-fold symmetry
as the tile border and notice nothing can (or should)
be filled into the interior of the first iterate
of the tile \tileS{A++}{6},
see Figure~\ref{fig:d333333-terdragon-1}.
Note how the area drawing on the right leaves half
of the interior empty.
This is because the curve is just half of the terdragon.

\IMAGE{d333333-terdragon-2}
The CW (digon) tile covers the area of the terdragon,
see the middle of Figure~\ref{fig:d333333-terdragon-2}.
The image on the right
for the axiom \texttt{A++A++A!A--A--A}
has the same shape as the tile for the terdragon
(with axiom \tileS{F+}{3}, a lattice tile),
compare to the right of Figure~\ref{fig:d333333-terdragon-1}.

Such a bisection of curves can be done for
curves with 2-fold symmetry.
Another example would be the (only) curve of order 4
with map
\begin{verbatim}
  F |--> F+F00F-F
\end{verbatim}

%
%
%

\subsubsection*{A curve of order 61}

The map of the curve of order 61 is
{\scriptsize
\begin{verbatim}
  A |--> A++A!A++A!A++A++A!A++A!A++A++A!A++A!A++A++A!A++A++A!A++A!A++A++A!A+A++A++A!A++ \
         A!A++A++A!A++A++A--A-A0A0A!A++A!A++A++A!A++A!A++A++A!A++A!A++A++A!A0A0A0A0A
\end{verbatim}
}

\IMAGE{d333333-r61-var2-motif}

The motif of the curve is shown in
Figure~\ref{fig:d333333-r61-var2-motif},
its CCW tile in Figure~\ref{fig:d333333-r61-var2-tile-A3},
and its CW (digon) tile in Figure~\ref{fig:d333333-r61-var2-tile-digon}.

\IMAGE{d333333-r61-var2-tile-A3}

\IMAGE{d333333-r61-var2-tile-digon}


\clearpage
\subsubsection*{A coloring with two edge classes}\label{sect:d333333-AB}

\IMAGE{d333333-AB-grid}
The coloring (see Figure~\ref{fig:d333333-AB-grid})
is such that we have two CCW prototiles,
\tileS{A++}{3} and \tileS{B++}{3}.
The color is switched over all edges,
so the CW (digon) tile is \tileS{A!B}{1}.

We need two CCW tiles that fit together.
The borders of the tiles chosen are shown next to each other in
Figure~\ref{fig:d333333-AB-r13-var00-build}.
Filling them gives the curves.
The resulting curve-set has order 13 and maps
\begin{verbatim}
  A |--> A++A!B++B!A++A++A!B0B+A
  B |--> B++B!A++A!B++B!A++A+B++B++B!A+B!A0A+B
\end{verbatim}

\IMAGE{d333333-AB-r13-var00-build}
%

\IMAGE{d333333-AB-r13-var00-tile-minus-decomp}
\IMAGE{d333333-AB-r13-var00-tile-plus-decomp}
The self-similarity of the CW (digon) tile \tileS{A!B}{1}
is shown in Figure~\ref{fig:d333333-AB-r13-var00-tile-minus-decomp}.
The mutual self-similarity of the CCW tiles in
Figure~\ref{fig:d333333-AB-r13-var00-tile-plus-decomp}.

In Figure~\jjref{fig:d333333-terdragon-2} we have split the terdragon
into halves in the minimal coloring.
With the coloring here we can
split any curve from the (minimally colored)
triangle grid (with no double edges)
into halves.

\IMAGE{d333333-curve-r13-1-digon}

The digon tile \tileS{A!B!}{1} of such a curve of order 13
with maps
\begin{verbatim}
  A |--> A0A++A0A0A--A++A--A--A0A++A++A--A
  B |--> B++B--B--B0B++B++B--B++B0B0B--B0B
\end{verbatim}
is shown in
Figure~\ref{fig:d333333-curve-r13-1-digon}.
We have seen half of this tile (that is, curve \texttt{A}) in
Figure~\jjref{fig:d333333-curve-r13-1-solid}.
This is the splitting of the curve
shown in Figure~\jjref{fig:333333-curve-r13-1}.

The same splitting technique works for
the square grid with double edges and
the coloring into two colors giving a chessboard pattern.
We'll do the equivalent for the trihex grid next.


\FloatBarrier
\clearpage
\subsection{The trihexagonal grid}\label{sect:d3636}

\IMAGE{d3636-grid}

Here we use turns by 60 degrees ($2\pi/6$).

The minimal coloring of the grid has two edge classes,
see Figure~\ref{fig:d3636-grid}.

The CCW prototiles are
\tileS{F+}{6} (blue hexagons) and
\tileS{G++}{3} (red triangles),
The CW (digon) prototile is
\tileS{F!G}{1}.
The transitions are
\texttt{F+F}, \texttt{F0G}, \texttt{F--F}, \texttt{F!G},
\texttt{G++G}, \texttt{G0F}, \texttt{G-G}, and \texttt{G!F}.

\subsection*{A curve-set of order 13}

\IMAGE{d3636-r13-1-motif}
We start with the curve of order 13 on the trihex grid
without double edges with map
\begin{verbatim}
  F |--> F+F+F+F+F--F+F+F--F--F+F+F--F
\end{verbatim}
It is shown on the left in
Figure~\ref{fig:d3636-r13-1-motif}.

For the grid with double edges we need a curve \texttt{G}
that is the reversal of the curve \texttt{F}:
its map is obtained by reversing the map for \texttt{F}
and changing
all \texttt{+} into  \texttt{-} and
all \texttt{-} into  \texttt{+}:
\begin{verbatim}
  G |--> G++G-G-G++G++G-G-G++G-G-G-G-G
\end{verbatim}
Now \tileS{F!G}{1} is already the CW (digon) tile,
nothing needs to be filled in,
see the two images on the right in the figure.


\IMAGE{d3636-r13-digon-tiling}
Figure~\ref{fig:d3636-r13-digon-tiling}
shows the tiling of the plane by
(second iterates) of the digon tiles.
Alternatively, it shows the tiling of the plane
by the curve \texttt{F} we started with,
on the grid without double edges.

\FloatBarrier
\subsection*{A curve-set of order 9 with the Koch snowflake as CCW tiles}

\IMAGE{d3636-r9-koch}
We can create a curve-set such that the CCW tile \tileS{G++}{3}
has 6-fold symmetry.
The most simple example is the curve-set of order 9
with maps
\begin{verbatim}
  F |--> F+F+F!G++G++G!F--F+F+F
  G |--> G++G!F+F!G++G++G-G
\end{verbatim}
whose CCW tiles are shown in
Figure~\ref{fig:d3636-r9-koch}.
The shape of both tiles is that of the Koch snowflake,
even having reflection symmetry.

\FloatBarrier
\subsection*{Another curve-set with two 6-fold symmetric tiles}

\IMAGE{d3636-rx-6fold-it1}

The borders of the CCW tiles of a curve-set
are shown in Figure~\ref{fig:d3636-rx-6fold-it1}.
The maps for the borders are
\begin{verbatim}
   F |--> F+F+F+F--F--F+F+F+F--F
   G |--> G++G-G-G-G++G++G-G-G-G
\end{verbatim}

\IMAGE{d3636-rx-6fold-it}
Again both tiles have 6-fold symmetry.
They do not have the reflection symmetry like the Koch snowflake, though.
This can bee seen in Figure~\ref{fig:d3636-rx-6fold-it}.


\FloatBarrier
\subsection{The \texorpdfstring{$(4.8.8)$}{(4.8.8)}-grid}\label{sect:d488-x}

Here we use turns by 45 degrees ($2\pi/8$).

\IMAGE{d488-x-grid-graph}
We turn all edges of the $(4.8.8)$-grid
into double edges as shown in the left in
Figure~\ref{fig:d488-x-grid-graph}.
There are three edge classes.
Letters \texttt{A} lie between octagons,
letters \texttt{b} and \texttt{B} lie
between a square and an octagon.
All edges are assigned an area left of them.
The areas for \texttt{b} fill the squares,
the areas for \texttt{A} and \texttt{B} fill the octagons,
see the right of Figure~\ref{fig:d488-x-grid-graph}.

The CCW prototiles are
\tileS{A+B+}{4} and \tileS{b++}{4},
CW prototiles are
\tileS{A!}{2} and \tileS{b!B!}{1}.
The transitions are
\texttt{A+B}, \texttt{A!A}, \texttt{A-b},
\texttt{b++b}, \texttt{b-A}, \texttt{b!B},
\texttt{B+A}, \texttt{B--B}, and \texttt{B!b}.

\IMAGE{d488-x-grid-graph-2}
Figure~\ref{fig:d488-x-grid-graph-2}
shows an area drawing of the grid where the colors
correspond to the orientation of the edges.


\subsubsection*{A curve of order 5 with two constant letters}

\IMAGE{d488-x-const-b-r5-var1}
The curve has order 5 and maps
\begin{verbatim}
    A |--> A-b++b++b++b!B+A+B--B+A!A+B+A
    b |--> b
    B |--> B
\end{verbatim}
Its first, second, and third iterates are shown in
Figure~\ref{fig:d488-x-const-b-r5-var1}.
The coloring used in the area drawing is
by orientation of the edges, as in
Figure~\ref{fig:d488-x-grid-graph-2}.
Note how all but the green and blue areas
(corresponding to the letter \texttt{A})
vanish with higher iterates.

\IMAGE{d488-x-const-b-r5-var1-tile-minus}
Iterates of the CW (digon) tile are shown in
Figure~\ref{fig:d488-x-const-b-r5-var1-tile-minus}.
It looks suspiciously like what is shown in
Figure~\jjref{fig:r5-dragon}.
Indeed, by dropping letter \texttt{b} and \texttt{B}
and normalizing turns one obtains the map
\begin{verbatim}
    A |--> A++A0A!A++A
\end{verbatim}
This is the map for the curve of order 5 on the
double-edged square grid from Section~\jjref{sect:d4444}.
Not a surprise,
dropping edges \texttt{b} and \texttt{B}
in our grid gives the square grid with double edges.


\FloatBarrier
\subsubsection*{A curve of order 13}

\IMAGE{d488-x-r13-2-var0-motifs}

The first iterates for a curve-set of order 13
with maps
\begin{verbatim}
  A |--> A+B!b++b++b++b-A+B+A!A+B+A+B!b-A+B+A
  b |--> b++b++b!B+A-b
  B |--> B+A+B+A+B!b++b++b-A+B!b-A+B--B+A+B
\end{verbatim}
are shown in Figure~\ref{fig:d488-x-r13-2-var0-motifs}.
Start and end points are marked.
Figure~\ref{fig:d488-x-r13-2-var0-curve-tiling}
shows how the curves tile the plane.

\IMAGE{d488-x-r13-2-var0-curve-tiling}




\FloatBarrier
\subsection{The \texorpdfstring{$(3.12.12)$}{(3.12.12)}-grid}
%
%
\IMAGE{d31212-grid}
The grid is shown in Figure~\ref{fig:d31212-grid}.
There are three edge classes.

Here we use turns by 30 degrees ($2\pi/12$).
The prototiles are
\tileS{A+B+}{6} and \tileS{b++++}{3} (CCW),
\tileS{A!}{2} and \tileS{b!B!}{1} (CW).
The transitions are
\texttt{A+B},
\texttt{A-b},
\texttt{A!A},
\texttt{B+A},
\texttt{B----B},
\texttt{B!b},
\texttt{b++++b},
\texttt{b-A}, and
\texttt{b!B}.
The letters at start and end of the transitions
determine the turn uniquely, so we can omit them
in the maps.

Our curve-set is of order 27 and has maps
\begin{verbatim}
  A |--> AB!bbbAB!bbbA!ABA!ABABABA!AB!bbbABABA
  B |--> BABABB!bbbABABB!bbbABB!bbbA!ABABABB!bbbABABAB
  b |--> bb!BAbbb!BAb
\end{verbatim}
We left out the turns except for the U-turns
\texttt{A!A}, \texttt{B!b}, and \texttt{b!B}.

\IMAGE{d31212-r27-var0-curves}
\IMAGE{d31212-r27-var0-curves-it}

The Motifs and third iterates are respectively shown in
Figures~\ref{fig:d31212-r27-var0-curves} and
\ref{fig:d31212-r27-var0-curves-it}.

\IMAGE{d31212-r27-var0-digons}

The CW (digon) tiles are shown in
Figure~\ref{fig:d31212-r27-var0-digons}.

\IMAGE{d31212-r27-var0-tiles-plus-small}
The CCW tiles, decomposed into curves, are shown in
Figure~\ref{fig:d31212-r27-var0-tiles-plus-small}.
Note that the CCW tile \tileS{b++++}{3}
has 6-fold symmetry.

\FloatBarrier

\IMAGE{d31212-r27-var0-tiles-plus-big-b3}
\IMAGE{d31212-r27-var0-tiles-plus-big-AB6}

Figures~\ref{fig:d31212-r27-var0-tiles-plus-big-b3}
and
\ref{fig:d31212-r27-var0-tiles-plus-big-AB6}
respectively show the tiles
\tileS{b++++}{3} and \tileS{A+B+}{6},
decomposed into smaller copies of each other.
Coloring of the small copies is just as in
Figure~\ref{fig:d31212-r27-var0-tiles-plus-small}.
The tiling of the plane by either the CCW tiles
or the CW (digon) tiles
(as in Figure~\ref{fig:d31212-r27-var0-digons})
can also be gleaned from
Figure~\ref{fig:d31212-r27-var0-tiles-plus-big-AB6}.


%
%
%


\FloatBarrier
\clearpage

\section{Curve-sets not fitting into our framework}

Curve-sets as we defined them include all edge-covering families
of plane-filling curves on grids we have seen before.
Still, there is more.
We present examples of edge-covering curves
on grids that do not fit into our framework.

In the following we refer to the objects presented as curve-sets
even though they are not curve-sets in the sense used so far.

\subsection{First example}

\IMAGE{nofit-manta-triangle}
On the triangle grid consider the curve-set given by the maps
\begin{verbatim}
  F |--> F+F-F-F0F+F+F-F0F
  H |--> H+F-F-F+H0H
\end{verbatim}
Turns are by 120 degrees ($2\,\pi/3$),
these curves live on the triangle grid.
Motif and second iterates are shown in Figure~\ref{fig:nofit-manta-triangle}.
Note that there are 12 letters \texttt{F} and 3 letters \texttt{H}
in the productions, so our concept of an order of a curve-set
does not apply here.
Edges \texttt{H} appear only on a line.

\IMAGE{nofit-manta-triangle-tiles}

Note that the transitions are not
as we required them for our curve-sets:
in the second iterate for \texttt{H} most of the transitions are
\texttt{F+F} and \texttt{F-F} (from the map for \texttt{F})
but some are \texttt{F-H} (from the map for \texttt{H}).
Therefore curve \texttt{H} does not live on a grid coloring
as the one we used.

%
The tile \tileS{F+}{3} is shown on the left in
Figure~\ref{fig:nofit-manta-triangle-tiles}
(the tile \tileS{F-}{3} is identical, modulo a reflection).
Curve \texttt{F} occupies \emph{half} of the plane,
its tile \tileS{H+}{3} is shown on the right.
The tile \tileS{H-}{3} has an empty triangular interior,
it is obtained by reflecting the three curves
in the tile \tileS{H+}{3}
over their straight borders.

\FloatBarrier
\subsection{Second example}

This is actually an infinite family of examples.

\IMAGE{nofit-r07-ex00-shape}
We start with a curve of order 7 on the triangle grid
given by the map
\begin{verbatim}
  F |--> F+F-F-F+F+F-F
\end{verbatim}
Its motif and decomposition into 7 smaller copies of itself
are shown in
Figure~\ref{fig:nofit-r07-ex00-shape}.

Now duplicate the map for \texttt{F} into one for the letter \texttt{G}.
Then, in the production,
replace the letter \texttt{F} by \texttt{G} at will.
We use the maps
\begin{verbatim}
  F |--> F+F-F-G+F+F-F
  G |--> G+G-G-G+G+G-G
\end{verbatim}
We have 6 letters \texttt{F} and 8 letters \texttt{G} in the production.
So this curve-set does not fit into our framework.

\IMAGE{nofit-r07-ex00-F}
These two curves \eqq{work} because they are of the same shape.
The production of \texttt{G} does not contain any letter \texttt{F}
so the curve is exactly the thing we started with.
Curve \texttt{F} looks more interesting, see Figure~\ref{fig:nofit-r07-ex00-F}:
the grid coloring induced by this curve has definitely no translation symmetry.

\IMAGE{nofit-r07-ex01-F}
We can do more by using
for \texttt{F} a curve of the same shape as \texttt{G}
(but not the identical sequence of turns):
\begin{verbatim}
  F |--> F0F+F+G-F-F0F  (changed)
  G |--> G+G-G-G+G+G-G
\end{verbatim}
Compare the motifs (top left) of Figure~\ref{fig:nofit-r07-ex01-F}
with that of Figure~\ref{fig:nofit-r07-ex00-F}!
In our illustration we use rounded corners to make regions with
many turns (red) visually distinguishable from regions with
straight lines.



\FloatBarrier

\vspace*{30mm}
\section*{Acknowledgments}

It is our pleasure to thank the following people
for their support, criticism, corrections, and improvements.
Michel Dekking,
Marvin Dittrich,
Daniel Fischer,
Christoph Haenel,
Edith Parzefall,
Neil Sloane,
and
Jeffrey Ventrella.
%
%


\clearpage
\makeatletter
\renewenvironment{theindex}%
{\columnseprule \z@
 \columnsep 21\p@
 \begin{multicols}{2}[\section*{Index of terms}%
   \addcontentsline{toc}{section}{\protect{Index of terms}}%
   \@mkboth{\indexname}%
   {\indexname}]%
    \thispagestyle{plain}\parindent\z@
    \parskip\z@ \@plus .3\p@\relax
    \let\item\@idxitem
    }%
{%
\end{multicols}\clearpage}
\makeatother

{
\baselineskip1.5em
\printindex
}
%
%

%
\FloatBarrier
\clearpage

\clearpage

\newcommand{\jjbibtitle}[1]{{\small\bfseries #1}}
\newcommand{\bdate}[1]{(#1)}


\end{document}

